# A Riemann-Roch theorem for flat bundles, with values in the algebraic Chern-Simons theory

By Spencer Bloch and Hélène Esnault*

### Introduction

Our purpose in this paper is to continue the algebraic study of complex local systems on complex algebraic varieties. We prove a Riemann-Roch theorem for these objects using algebraic Chern-Simons characteristic classes.

A complex local system $\mathcal{E}$ on a smooth, projective complex variety $X$ gives rise to a locally free analytic sheaf $E^{\mathrm{an}} := \mathcal{E} \otimes_{\mathbb{C}} \mathcal{O}_X^{\mathrm{an}}$ which (using GAGA) admits a canonical algebraic structure $E$. The tautological analytic connection on $\mathcal{E} \otimes_{\mathbb{C}} \mathcal{O}_X^{\mathrm{an}}$ induces an integrable algebraic connection $\nabla : E \to E \otimes \Omega_X^1$. Combining GAGA with the Poincaré lemma, we see that the analytic cohomology of the local system can be identified with the hypercohomology of the algebraic de Rham complex

$$(0.1) \qquad \Omega_X^* \otimes_{\mathcal{O}_X} E := \{ E \xrightarrow{\nabla} E \otimes \Omega_X^1 \xrightarrow{\nabla} E \otimes \Omega_X^2 \xrightarrow{\nabla} \cdots \}.$$

We will work with algebraic connections $\nabla : E \to E \otimes \Omega_X^1$ where $X$ is an algebraic variety defined over a field $k$ of characteristic 0. To understand what kind of Riemann-Roch theorem we might expect for such objects, we may apply the Grothendieck-Riemann-Roch, with Chern classes in the Chow group modulo torsion to a relative version of (0.1). This sort of calculation was first done by Mumford [23]. Using the remarkable identity ([19]) for a vector bundle $T$ of rank $d$,

$$(0.2) \qquad (-1)^d c_d(T^*) = \mathrm{Td}(T) \cdot \sum_i (-1)^i \mathrm{ch}(\wedge^i T^*),$$

one deduces for $f : X \to S$ smooth and proper of fibre dimension $d$, that

$$(0.3) \qquad \mathrm{ch}(\mathbb{R}f_*(E \otimes \Omega_{X/S}^*)) = (-1)^d f_*(\mathrm{ch}(E) \cdot c_d(\Omega_{X/S}^1)).$$

The situation is not totally satisfactory because the above Riemann-Roch depends only on the graded sheaf $E \otimes \Omega^*$ and does not depend on the connection $\nabla$. A theory of algebraic differential characters $\mathrm{AD}(X)$ and characteristic classes $c_p(E, \nabla) \in \mathrm{AD}^p(X)$ was developed in [12], [14]. These classes simultaneously refine the Chow group classes $c_p(E)$ and the differential character classes of Chern and Simons [7]. It seems likely that the optimal Riemann-Roch theorem for flat bundles (and perhaps more generally for regular holonomic $\mathcal{D}$-modules) will take values in this theory. In the present work, characteristic classes $w_p(E, \nabla)$ (analogous to Chern classes) and $\mathrm{N}w_p(E, \nabla)$ (Newton classes, analogous to the Chern character) lie in a quotient of $\mathrm{AD}(X)$ which we call Chern-Simons cohomology:

$$(0.4) \qquad \mathrm{N}w_1(E, \nabla) \quad = \quad \mathrm{N}w_1(\det(E), \nabla) \in H^2_{\mathrm{CS}}(X)$$

$$(0.5) \qquad \qquad \qquad := \quad H^0(X, \Omega^1_X / d\log(\mathcal{O}^*_X))$$

$$\mathrm{N}w_p(E, \nabla) \in H^{2p}_{\mathrm{CS}}(X) \quad := \quad H^0(X, \Omega^{2p-1}_X / d\Omega^{2p-2}_X), \ p \geq 2.$$

---
*This work has been partly supported by NSF grant DMS-9423007-A1 and the DFG Forschergruppe "Arithmetik und Geometrie."

One has

$$H_{\mathrm{CS}}(X) \hookrightarrow H_{\mathrm{CS}}(k(X)) \cong \mathrm{AD}(k(X)),$$

so that algebraic Chern-Simons characteristic classes can be thought of as those parts of the AD-classes which survive at the generic point. A connection on $\mathcal{O}_X$ is determined by a 1-form $\eta = \nabla(1)$, and

$$w_1(E, \nabla) = \mathrm{N}w_1(\mathcal{O}_X, \nabla) := \eta \bmod d\log(\mathcal{O}^*_X).$$

In general, $\mathrm{N}w_1(E, \nabla) = 0$ if and only if $\det(E)$ has a nontrivial flat section.

The $\mathrm{N}w_p$ for $p \geq 2$ are related to the classes $w_p$ described in [5], (0.2.3), in the same way that the Chern character is related to the Chern class. Zariski-locally the bundle $E$ is trivial. Write $A$ for a locally defined connection matrix, and let $F(tA) = tdA - t^2 A \wedge A$ be the curvature of the connection $tA$. Define

$$(0.6) \qquad \qquad TP(A) = p \int_0^1 P(A \wedge F(tA)^{p-1}) dt$$

where $P$ is an invariant polynomial of degree $p$ on the Lie algebra. One has

$$(0.7) \qquad \qquad dTP(A) = P(F(A)).$$

A gauge transformation $A \mapsto gAg^{-1} + dgg^{-1}$ changes $TP(A)$ by a Zariski-locally exact form, so that these forms glue to a section

$$(0.8) \qquad w(E, P, \nabla) \in H^{2p}_{\mathrm{CS}}(X) := H^0(X, \Omega^{2p-1}_X / d\Omega^{2p-2}_X).$$

The classes $\mathrm{N}w_p$ are obtained by taking $P(M) = \mathrm{Tr}(M^p)$. When we think of the Chern character, $\mathrm{N}w_p$ appears to be more natural in a Riemann-Roch problem than the classes $w_p$ obtained by taking $P$ to be the invariant polynomial of degree $p$ whose value on diagonal matrices is the $p^{\mathrm{th}}$ elementary symmetric function. The N in the terminology stands for Newton. When the connection is integrable it follows from (0.7) that $w$ is closed; i.e.,

$$(0.9) \qquad w(E, P, \nabla) \in H^0(X, \mathcal{H}^{2p-1}) = \ker(H^{2p}_{\mathrm{CS}}(X) \xrightarrow{d} H^0(X, \Omega^{2p}_X)).$$

Here $\mathcal{H}^r$ is the $r^{\mathrm{th}}$ cohomology sheaf of the algebraic de Rham complex $\Omega^*_X$.

Globally defined connections, even nonflat ones, are rather rare. For applications, it is important to work with connections admitting log poles along a normal crossings divisor $Y \subset X$:

$$(0.10) \qquad \nabla : E \to E \otimes \Omega^1_X(\log Y).$$

By definition [8], sections of $\Omega^1_X(\log Y)$ are meromorphic forms $\omega$ such that $f\omega$ and $f\,d\omega$ are regular, where $f$ is a local defining equation for $Y$. One defines (Definition 1.8 below) classes

$$\mathrm{Nw}_p(E, \nabla) \in H^{2p}_{\mathrm{CS}}(X(\log Y)).$$

Let $f : X \to S$ be a flat map of smooth varieties, and let $Y \subset X$ and $T \subset S$ be normal crossings divisors. The data

$$\{f : X \to S, Y \subset X, T \subset S\}$$

are said to be relative normal crossings divisors if $f^{-1}(T) \subset Y$ and

$$\Omega^1_{X/S}(\log(Y)) := \Omega^1_X(\log(Y))/f^*\Omega^1_S(\log(T))$$

is locally free of rank $d = \dim(X/S)$. Let $\{Z_i\}$ be the components of $Y$ not lying in $f^{-1}(T)$. The residue map

$$\mathrm{res}_{Z_i} : \Omega^1_{X/S}(\log(Y)) \to \mathcal{O}_{Z_i}$$

is surjective; thus in particular the top Chern class $c_d(\Omega^1_{X/S}(\log(Y))|_{Z_i}$ vanishes in the corresponding Chow group $CH^d(Z_i)$. We call $\mathrm{res}_{Z_i}$ a partial trivialization of $\Omega^1_{X/S}(\log(Y))$. Using ideas of T. Saito [27], [28], one can define relative top Chern classes in a relative Chow group

$$(0.11) \quad c_d(\Omega^1_{X/S}(\log(Y)), \mathrm{res}_Z) \;\in\; CH^d(X, Z_\bullet)$$
$$:= \; \mathbb{H}^d(X, \underline{K}_{d,X} \to \underline{K}_{d,Z^{(1)}} \to \underline{K}_{d,Z^{(2)}} \to \dots).$$

(Here $\underline{K}_d$ denotes the $d^{\mathrm{th}}$ Milnor $K$-sheaf, defined (Definition 1.6) to be the image of the sheaf of symbols in the Milnor $K$-theory of the function field. Alternatively, one can interpret $\underline{K}_d$ as the Zariski sheaf associated to the higher Chow group $CH^d(X, d)$ ([4]) and $Z^{(i)}$ is the normalized $i$-fold intersection of components of the normal crossings divisor $Z = Y - f^{-1}(T)$.) One has pairings and a trace map

$$(0.12) \quad H^{2p}_{\mathrm{CS}}(X(\log Y)) \times CH^q(X, Z_\bullet) \;\to\; H^{2p+2q}_{\mathrm{CS}}(X(\log f^{-1}(T)))$$
$$\xrightarrow{f_*} \; H^{2p+2q-2d}_{\mathrm{CS}}(S(\log(T))).$$

The final ingredient needed to formulate the Riemann-Roch theorem is the existence of a canonical Gauss-Manin connection on the de Rham cohomology,

$$(0.13) \qquad \nabla_{\mathrm{GM}} : H^i_{\mathrm{DR}}(X/S(\log Y), E) := \mathbb{R}^i f_*(E \otimes \Omega^*_{X/S}(\log(Y)))$$
$$\to \; H^i_{\mathrm{DR}}(X/S(\log Y), E) \otimes \Omega^1_S(\log T).$$

Traditionally, $\nabla_{\mathrm{GM}}$ is defined under the hypothesis that $\nabla$ on $E$ is integrable, but in fact one need only assume that the curvature form is basic (Definition 3.1); i.e.,

$$(0.14) \qquad \nabla^2 \in \mathrm{Hom}(E, f^*\Omega_S^2(\log T) \otimes E).$$

THEOREM 0.1 (Riemann-Roch). *Let $\{f : X \to S, Y \subset X, T \subset S\}$ be a relative normal crossings divisor with $f$ flat and projective and $\dim(X/S) = d$. Let $(E, \nabla)$ be a locally free sheaf with connection on $X$, and assume the curvature form $\nabla^2$ is basic (0.14). Then*

$$(0.15) \quad \mathrm{N}w_p(H^*_{\mathrm{DR}}(X/S(\log Y), E), \nabla_{\mathrm{GM}})$$
$$- \mathrm{rank}(E) \cdot \mathrm{N}w_p(H^*_{\mathrm{DR}}(X/S(\log Y)), \nabla_{\mathrm{GM}})$$
$$= (-1)^d f_*\left( c_d(\Omega^1_{X/S}(\log Y), \mathrm{res}_Z) \cdot \mathrm{N}w_p(E, \nabla) \right).$$

*Here $H^*_{\mathrm{DR}}(X/S(\log Y))$ is the usual de Rham cohomology ($E = \mathcal{O}_X$ with the trivial connection). Now,*

$$\mathrm{N}w_p(H^*_{\mathrm{DR}}(X/S(\log Y)), \nabla_{\mathrm{GM}}) = (0), \ p \geq 2,$$
$$2 \cdot \mathrm{N}w_1(H^*_{\mathrm{DR}}(X/S(\log Y)), \nabla_{\mathrm{GM}}) = (0).$$

*Remark* 0.2. When $\nabla$ is flat, the theorem remains true with $\mathrm{N}w$ replaced by $w$ (Corollary 1.10).

*Remark* 0.3. Note that this is really a Riemann-Roch for virtual bundles of rank zero. It would be of interest to give a "Noether Formula" describing

$$(0.16) \quad \mathrm{N}w_1(H^*_{\mathrm{DR}}(X/S(\log Y)), \nabla_{\mathrm{GM}}) \in H^2_{\mathrm{CS}}(S(\log T))_2 \cong H^1_{\text{ét}}(S - T, \mathbb{Z}/2\mathbb{Z})$$

in terms of characteristic classes for $X/S$. More precisely, one would like to understand the residue of this class along a component of $T$ in terms of suitable characteristic classes with support, supported over the component.

Working analytically with the local system $E^\nabla$ (and $Y = \emptyset$), Bismut-Lott [3] and Bismut [2] proved an analogue of Theorem 0.1 using characteristic classes $\hat{c}_n(E^\nabla) \in H^{2n-1}(X_{\mathrm{an}}, \mathbb{C}/\mathbb{Q}(n))$ as defined by Chern and Simons [7].

Analogously, if $S$ is replaced by a finite field $\mathbb{F}_q$, $(E, \nabla)$ by a tame representation of the fundamental group $\rho$, then Deligne's Theorem [9], [10], and subsequent work by Laumon [21], S. Saito [26] and T. Saito [27] show that for $n = 1$ the formula (0.15) remains true

$$(0.17) \qquad \bigotimes_i (\det H^i_{\text{ét}}(X, \rho - \dim(\rho) \cdot \mathbb{Q}_\ell))^{(-1)^i} = (\det \rho | c_d(\Omega^1_X))^{(-1)^d}$$

as dimension 1, $\mathbb{Q}_\ell$-vectorspaces with Frobenius action. Deligne's proof [10] in the case $d = 1, S = \mathbb{F}_q$ and $\rho$ is a character of the fundamental group is purely geometric, and relies on properties of $\mathrm{Pic}^0(X)$. In view of the shape of the formula (0.15) involving only $c_d(\Omega^1_{X/S})$, it would be natural to try to use Higgs cohomology and the geometry of the Hitchin map to recover Theorem 0.1. We could not do this. On the other hand, $\mathrm{GM}^i(\nabla)$ is a particular case of the image of a $\mathcal{D}$-module by a projective morphism. Thus Theorem 0.1 should have a formulation for images of regular holonomic $\mathcal{D}$-modules under projective maps. To this end, an algebraic Chern-Simons and Cheeger-Simons theory of characteristic classes of such $\mathcal{D}$-modules remains to be constructed. Note that a Riemann-Roch theorem for the category of all (not necessarily holonomic) $\mathcal{D}$-modules is known [22]. It has a completely different flavor.

It should be stressed that, although their definition is simple and natural, the higher algebraic Chern-Simons classes $\mathrm{N}w_p(E)$, $p \geq 2$, are mysterious. Working with bundles with basic curvature forms, we give examples where these classes are nonzero, so that our Riemann-Roch theorem has some content in degrees $\geq 2$. However, we hope the methods developed in this paper can be extended to yield a Riemann-Roch theorem with values in the algebraic differential characters

(0.18)
$$\mathrm{AD}^n(X(\log Y)) := \mathbb{H}^n\left(X, \underline{\underline{K}}_n \overset{d\log}{\to} \Omega^n_X(\log Y) \to \ldots \to \Omega^{2n-1}_X(\log Y)\right),$$

refining the theorem given here and also the analytic work of [3] and [2].

Briefly, Section 1 establishes some necessary technical results about a Grothendieck group of modules with connections and introduces the Chern-Simons groups where our characteristic classes take values. Section 2 discusses relative top Chern classes, products, and covariant functoriality. Arguments in Section 3 reduce the Riemann-Roch theorem to the case of $\mathbb{P}^1$ (with log poles), and Section 4 establishes the theorem in that case. Finally, in Section 5, we prove the Riemann-Roch theorem for curves over function fields and the class $\mathrm{N}w_1$ (=determinant) without the basic curvature hypothesis.

We are endebted to A. Beilinson and T. Saito for considerable help and to O. Gabber who pointed out an error in an early version of the paper and provided us with some unpublished notes of his own on the theorem of Bolibruch [17].

## 1. Grothendieck group of coherent sheaves with connections

Our purpose in this section is to make sense of the expression

$$\mathrm{N}w_p(\mathbb{R}f_*(E \otimes \Omega^*_{X/S}(\log D)), \nabla_{\mathrm{GM}})$$

which appears in the Riemann-Roch theorem. Rather than develop the notion of objects in the derived category with connections and their characteristic classes, it is simpler to take

$$\sum_i (-1)^i \mathrm{N}w_p(\mathbb{R}^i f_*(E \otimes \Omega^*_{X/S}(\log D)), \nabla_{\mathrm{GM}}).$$

The difficulty with this is that, because our connections have log poles, the coherent sheaves appearing on the left need not be locally free. We will define a Grothendieck group of sheaves with connections which is large enough to contain expressions like the above, and on which the $\mathrm{N}w_p$ are defined. Since in fact, one can even define classes of connections in the group of algebraic, differential characters, after enlarging the divisor of poles, we insert the construction in this section. For the main Theorem 0.1, we need only Proposition 1.5 and Corollary 1.11 ii), and not the whole strength of Proposition 1.4 and Corollary 1.11 i).

LEMMA 1.1.   *Let $X$ be a smooth variety over a function field $F$ over an algebraically closed field $k$ of characteristic 0; let $D \subset X$ be a normal crossings divisor,   and let $\Omega^1_X(\log D) := \Omega^1_{X/k}(\log D)$ be the locally free sheaf of differential forms with logarithmic poles along $D$ (see Definition 2.1).  Then there is a normal crossings divisor $Y \subset X$ such that*

   (i) *$Y = D + H$ for some normal crossings divisor $H$,*

   (ii) *all irreducible components of $H$ are very ample,*

   (iii) *$\Omega^1_X(\log Y)$ is generated by global sections.*

*Proof.* Fix a line bundle $L$ on $X$ such that $L$ and $L(-D_i)$ are very ample for all $i$. We take $H = \sum_j H_j$ where the $H_j$ are defined by sections of $L$ and $L(-D_i)$. We suppose that $Y = H + D$ is a normal crossings divisor, and the $\{H_j\}$ contain all the coordinate hyperplanes for some general set of projective coordinates on $\mathbb{P}(\Gamma(X, L))$.  Global forms $df/f$ for $(f) = H_j - H_\ell$ suffice to generate $\Omega^1_X \subset \Omega^1_X(\log Y)$.  Assume further that for any $i$ and any $x \in D_i$, there exists an $H_j$ defined by a section of $L(-D_i)$ such that $x \notin H_j$. Take $H_\ell$ in the linear system defined by $L$ such that $x \notin H_\ell$. Then $D_i + H_j - H_\ell = (f)$ for some rational function $f$, and $df/f \in H^0(X, \Omega^1_X(\log Y))$ has residue 1 along $D_i$ and no other residue through $x$.                                             □

*Definition* 1.2.   Let $D \subset X$ be a normal crossings divisor on a smooth variety defined over a field $F$.  Let $\Omega := \Omega^1_{X/k}(\log D)$ for some $k \subset F$.  A connection $\nabla$ on a coherent $\mathcal{O}_X$-module $M$ is a $k$-linear map

$$\nabla : M \to M \otimes \Omega$$

satisfying the Leibniz rule $\nabla(xm) = x\nabla(m) + m \otimes dx$. A short exact sequence of connections is a commutative diagram

(1.1)
$$
\begin{array}{ccccccccc}
0 & \longrightarrow & M' & \longrightarrow & M & \longrightarrow & M'' & \longrightarrow & 0 \\
& & \nabla' \downarrow & & \nabla \downarrow & & \nabla'' \downarrow & & \\
0 & \longrightarrow & M' \otimes \Omega & \longrightarrow & M \otimes \Omega & \longrightarrow & M'' \otimes \Omega & \longrightarrow & 0
\end{array}
$$

of exact sequences. The curvature $\nabla^2 : M \to \Omega^2_X(\log D) \otimes M$ is the $\mathcal{O}_X$-linear map induced by $\nabla \circ \nabla$ and the projection $\Omega \otimes \Omega \to \Omega^2_X(\log D)$. We do not assume $M$ is locally free.

*Definition* 1.3. (i) Define $\mathcal{G}(X, \log D)$ to be the Grothendieck group of coherent sheaves with connections on $X$ as above, with relations

$$[M, \nabla] = [M', \nabla'] + [M'', \nabla'']$$

coming from short exact sequences of modules with connection.

(ii) Define $\mathcal{K}(X, \log D)$ to be the corresponding object, where modules $M$ are required to be locally free.

PROPOSITION 1.4. *Let $X$ be projective and smooth over a function field $F$. Assume that $\Omega = \Omega_X(\log Y)$ is generated by global sections and that $Y$ contains an ample irreducible component. Then the natural map $\mathcal{K}(X, \log Y) \to \mathcal{G}(X, \log Y)$ is an isomorphism.*

*Proof.* Let $H \subset Y$ be an ample irreducible component. Then for any $p \in \mathbb{Z}$, the sheaf $\mathcal{O}(pH)$ carries a canonical connection with residue $-p \cdot$Identity along $H$ induced by the trivial connection $(\mathcal{O}_X, d)$. For $(M, \nabla)$ a connection, one has then a tensor connection on $M(pH) := M \otimes \mathcal{O}_X(pH)$, still denoted by $\nabla$, with residue $\mathrm{Res}_H(\nabla) - p \cdot$Identity along $H$. Choose $r \geq 0$ sufficiently large such that $M(rH)$ is generated by global sections. Define a locally free sheaf $A$ by the exact sequence

$$0 \to A \to \Gamma(X, \Omega) \otimes \mathcal{O}_X \to \Omega \to 0.$$

We may further suppose $H^1(X, A \otimes M(rH)) = (0)$. Tensoring the above sequence with $M(rH)$ and taking global sections, we see that

$$\Gamma(X, \Omega) \otimes \Gamma(X, M(rH)) \twoheadrightarrow \Gamma(\Omega \otimes M(rH)).$$

Let $e_i$ be a basis of $\Gamma(X, M(rH))$, and choose

$$\sum \omega_{ij} \otimes e_j \in \Gamma(X, \Omega) \otimes \Gamma(X, M(rH))$$

lifting $\nabla(e_i)$. Define a connection $\Phi$ on $\mathcal{O}_X \otimes \Gamma(X, M(rH))$ by

$$\Phi(1 \otimes e_i) = \sum \omega_{ij} \otimes e_j.$$

The commutative diagram (defining $K$ and $\Psi$)

$$(1.2)$$

$$
\begin{array}{ccccccccc}
0 & \longrightarrow & K(r) & \longrightarrow & \mathcal{O}_X \otimes \Gamma(X, M(r)) & \longrightarrow & M(r) & \longrightarrow & 0 \\
& & {\scriptstyle V\Psi}\downarrow & & {\scriptstyle \Phi}\downarrow & & {\scriptstyle \nabla}\downarrow & & \\
0 & \longrightarrow & \Omega \otimes K(r) & \longrightarrow & \Omega \otimes \Gamma(X, M(r)) & \longrightarrow & \Omega \otimes M(r) & \longrightarrow & 0
\end{array}
$$

can be twisted by $\mathcal{O}_X(-rH)$ (which does not destroy the log connection) to yield the first step in a resolution of $M$ by locally free sheaves with connection. By induction on the homological dimension of $M$, we conclude that such resolutions exist.

Given two such resolutions $P^{\cdot}$ and $Q^{\cdot}$ of $M$, we must show

$$\sum (-1)^i [P^i, \nabla_{P^i}] = \sum (-1)^i [Q^i, \nabla_{Q^i}] \in \mathcal{K}(X, \log D).$$

If we have a surjection $Q^{\cdot} \twoheadrightarrow P^{\cdot}$ compatible with connections, this is clear since the kernel complex is locally free and acyclic and so represents 0 in $\mathcal{K}(X, \log D)$. It therefore suffices to show in general we can construct a resolution $R^{\cdot}, \nabla_R$ of $M$, and surjections $R \twoheadrightarrow P^{\cdot}, Q^{\cdot}$. We choose as above

$$(R^0, \nabla_{R^0}) \twoheadrightarrow (P^0 \amalg_M Q^0, \nabla_{P^0 \amalg Q^0}).$$

Suppose we have constructed $R^0, \cdots, R^{i-1}$. Now, there is a diagram

$$(1.3)$$

$$
\begin{array}{ccccccc}
0 & \longrightarrow & \ker(d_R) & \longrightarrow & R^{i-1} & \xrightarrow{d_R} & R^{i-2} \\
& & \downarrow & & \downarrow & & \downarrow \\
0 & \longrightarrow & \ker(d_P) & \longrightarrow & P^{i-1} & \xrightarrow{d_P} & P^{i-2}.
\end{array}
$$

Adding a summand to $R^{i-1}$ mapping to 0 in $R^{i-2}$ we may assume all three vertical maps are surjective. We have a similar diagram for $Q$, and again vertical maps can be taken to be onto. We construct $R^i$ with connection mapping onto the coproduct

$$P^i \amalg_{\ker(d_P)} \ker(d_R) \amalg_{\ker(d_Q)} Q^i.$$

This coproduct surjects onto $P^i, Q^i$, and $\ker(d_R)$, so that the inductive step is complete.

It follows that the natural map $\rho : \mathcal{K}(X, \log Y) \to \mathcal{G}(X, \log Y)$ is surjective. A similar construction shows that an exact sequence of modules with connection can be lifted to an exact sequence of resolutions, so that $\rho$ admits a surjective splitting $[M, \nabla] \mapsto \sum (-1)^i [P^i, \nabla_{P^i}]$, and the two groups are isomorphic. $\quad\square$

PROPOSITION 1.5. *Let $X$ be a smooth affine variety over a field $k$, and $D \subset X$ be a normal crossings divisor. Then the natural map $\mathcal{K}(X, \log D) \to \mathcal{G}(X, \log D)$ is an isomorphism.*

*Proof.* The proof is similar to that of Proposition 1.4.                □

*Definition* 1.6.   If $F$ is a field, we let $K_n(F)$ be the group of Milnor $K$-theory of $F$. On a variety $X$, we let $\underline{K}_n$ be the image of the Zariski sheaf of Milnor $K$-theory into $\oplus_{x \in X^0} i_{x,*} K_n(k(x))$, where the $X^0$ are the points of $X$ of codimension 0 and $i_x : x \to X$. Alternatively, we may define $\underline{K}_n$ to be the sheaf associated to $U \mapsto CH^n(U, n)$ [4].

According to [24], [20], [4], [16], [25], the sheaves $\underline{K}_n$ on a smooth variety $X$ have all the good properties the Quillen $K$-sheaves have. Recall that the variety $X$ is always assumed to be smooth.

*Definition* 1.7.   We will also use the notation $\underline{\underline{G}}_{q,X}$ for the Gersten resolution of $\underline{\underline{K}}_{q,X}$.

Other good properties include the isomorphism $H^r_{\{x\}}(\underline{K}_n)$ with $K_{n-r}(k(x))$ where $x \in X^r$ is a point of codimension $r$, the projective bundle formula and the localization sequence. Moreover, there is a $d\log$ map $d\log : \underline{K}_n \to \Omega^n_X$, induced by $d\log : \underline{K}_1 \to \Omega^1_X$ by tensor product, since the kernel of the Milnor $K$-sheaf to $\oplus_{x \in X^0} i_{x,*} K_n(k(x))$ is supported on proper closed subsets, and thus is killed by differentiation.

One introduces now the following complex:

$$(1.4) \qquad \Omega\underline{K}_{n,X}(\log D) \quad := \quad j_* \underline{K}_{n,X-D} \xrightarrow{d\log} \Omega^n_X(\log D)$$
$$\to \Omega^{n+1}_X(\log D) \to \ldots \to \Omega^{2n-1}_X(\log D),$$

where $j : X - D \to X$ is the embedding. This complex differs slightly from

$$(1.5) \qquad \Omega'\underline{K}_{n,X}(\log D) \quad := \quad \underline{K}_{n,X} \xrightarrow{d\log} \Omega^n_X(\log D)$$
$$\to \Omega^{n+1}_X(\log D) \to \ldots \to \Omega^{2n-1}_X(\log D),$$

the complex used in [12].

One introduces as in [12] the groups

$$\mathrm{AD}^n(X, \log D) = \mathbb{H}^n(X, \Omega'\underline{K}_{n,X}(\log D)),$$

and as in [5] the corresponding Zariski sheaves

$$\underline{H}^n(\Omega\underline{K}_{n,X}(\log D)).$$

*Definition* 1.8.   The algebraic Chern-Simons groups are defined by

$$H^{2n}_{\mathrm{CS}}(X(\log D)) := \Gamma(X, \underline{H}^n(\Omega\underline{K}_{n,X}(\log D))).$$

PROPOSITION 1.9.   *The restriction maps to the generic point*

$$(1.6) \qquad H^2_{\mathrm{CS}}(X(\log D)) \quad \rightarrow \quad \Omega^1_{k(X)}/d\log k(X)^{\times}$$

$$H^{2n}_{\mathrm{CS}}(X(\log D)) \quad \rightarrow \quad \Omega^{2n-1}_{k(X)}/d\Omega^{2n-2}_{k(X)} \text{ for } n > 1$$

*are injective.*

*Proof.* Assume first $n > 1$. Let $R$ be the local ring at a point of $X$. It will suffice to show that the map

$$\Omega^{2n-1}_R(\log D)/d\Omega^{2n-2}_R(\log D) \rightarrow \Omega^{2n-1}_{k(X)}/d\Omega^{2n-2}_{k(X)}$$

is injective. Since the analogous inclusion on exact $2n$-forms is evident, one reduces to showing the injectivity of the map

$$H^{2n-1}_{\mathrm{DR}}(\mathrm{Spec}\,(R) - D) \rightarrow H^{2n-1}_{\mathrm{DR}}(k(X)).$$

This follows from the diagram

$$
\begin{array}{ccccc}
0 \longrightarrow H^{2n-1}_{\mathrm{DR}}(\mathrm{Spec}\,(R)) & \longrightarrow & H^{2n-1}_{\mathrm{DR}}(\mathrm{Spec}\,(R) - D) & \longrightarrow & \oplus_i H^{2n-2}_{\mathrm{DR}}(D_i) \\
\| & & \downarrow & & \downarrow \text{inject} \\
0 \longrightarrow H^{2n-1}_{\mathrm{DR}}(\mathrm{Spec}\,(R)) & \longrightarrow & H^{2n-1}_{\mathrm{DR}}(k(X)) & \longrightarrow & \oplus_i H^{2n-2}_{\mathrm{DR}}(k(D_i))
\end{array}
$$

which is a part of the Gersten style resolution for de Rham cohomology [6].

For $n = 1$, one observes that $d\log$ of a rational function has no poles along a divisor if and only if the rational function is regular invertible along this divisor. $\qquad \square$

COROLLARY 1.10.   *Let $X$ be a smooth variety over a field $k$, and $D \subset X$ be a normal crossings divisor. The functorial classes defined in [5] on $\mathcal{K}(X, \log D)$ extend to $\mathcal{G}(X, \log D)$:*

$$w_n(M, \nabla) \in H^{2n}_{\mathrm{CS}}(X(\log D)).$$

*Moreover, one defines "Newton" classes $\mathrm{N}w_n(M, \nabla) \in H^{2n}_{\mathrm{CS}}(X(\log D))$ by requiring that*

$$\mathrm{N}w_n(M, \nabla) = P_n(w_1(M, \nabla), \ldots, w_n(M, \nabla)) \in H^{2n}_{\mathrm{CS}}(X(\log D))$$

*where $P_n$ is the universal polynomial of degree $n$ with $\mathbb{Z}$-coefficients expressing the Newton classes in terms of the Chern classes (or, what is the same thing, expressing the symmetric function "sum of $n^{\mathrm{th}}$ powers" in terms of the elementary symmetric functions). In particular $w_1(M, \nabla) = \mathrm{N}w_1(M, \nabla)$.*

*Proof.* Let $X = \cup X_i$ be an affine covering of $X$. By Proposition 1.5, $w_n((M, \nabla)|X_i) \in H^0(X_i, \underline{H}^n(\Omega \underline{K}_{n,X}(\log D)))$ is well-defined on $X_i$, and verifies $w_n((M, \nabla)|X_i)|X_{ij} = w_n((M, \nabla)|X_j)|X_{ij}$ for $i \neq j$. Proposition 1.9 allows us to conclude. $\qquad \square$

As a corollary, one obtains the existence of well-defined algebraic differential characters for connections on coherent sheaves, after enlarging a bit the poles, and of well-defined Chern-Simons invariants, even without enlarging the poles.

COROLLARY 1.11. *Let $X$ be as in Lemma 1.1. Let $D$ be a normal crossings divisor, $Y = D + H$ be a normal crossings divisor such that $\Omega^1_X(\log Y)$ is globally generated. Let $(M, \nabla) \in \mathcal{G}(X, \log D)$. Denote by $\nabla_Y$ the same connection, but understood as a connection with logarithmic poles along $Y$. Then one has the following*:

(i) *The functorial and additive classes defined in* [12] *on $\mathcal{K}(X, \log D)$ extend to $\mathcal{G}(X, \log Y)$*:
$$c_n(M, \nabla_Y) \in \mathrm{AD}^n(X, \log Y).$$

(ii) *The functorial and additive classes defined in* [5] *on $\mathcal{K}(X, \log D)$ extend to $\mathcal{G}(X, \log D)$*:
$$w_n(M, \nabla) \in H^{2n}_{\mathrm{CS}}(X(\log D)) \subset H^{2n}_{\mathrm{CS}}(X(\log Y)),$$
*and in this larger group involving poles along $Y$, $w_n(M, \nabla)$ is the image of $c_n(M, \nabla_Y)$ under the natural map*
$$\mathrm{AD}^n(X, \log Y) \to H^{2n}_{\mathrm{CS}}(X(\log Y)).$$

*Proof.* By Propositions 1.4 and 1.9 and by [12], we must show that $w_n(M, \nabla_Y)$ has no poles along $H$. Let $Y' = D + H'$ be another normal crossings divisor such that $Y + H'$ is a normal crossings divisor, $H$ and $H'$ have no common component, and such that $\Omega^1_X(\log Y')$ is globally generated. Then

$$\begin{aligned} \mathrm{im}\, c_n(M, \nabla_Y) &= \mathrm{im}\, c_n(M, \nabla_{Y'}) \\ &= c_n(M, \nabla_{D+H+H'}) \in \mathrm{AD}^n(X, \log(D + H + H')). \end{aligned}$$

Thus

$$\begin{aligned} w_n(M, \nabla_Y) &= w_n(M, \nabla_{Y'}) \\ &\in H^{2n}_{\mathrm{CS}}(X(\log Y)) \subset H^{2n}_{\mathrm{CS}}(X(\log(Y + H + H'))) \end{aligned}$$

and therefore, $w_n(M, \nabla_Y)$ has no residues along $H$. $\qquad\square$

## 2. Relativity

In this section we consider various constructions involving relative normal crossings divisors, relative Chern classes, and related questions. In order to formulate our Riemann-Roch theorem, we need a notion of the top Chern class of the differential forms with logarithmic poles, together with a trivialization along the poles coming from the residue map.

Let $D = \sum_{i=1}^d D_i$ be a divisor on a smooth variety $X$. Assume the $D_i$ are irreducible, and write $D_I = \cap_{i \in I} D_i$. By convention, $D_\emptyset = X$.

*Definition* 2.1.   The divisor $D$ is said to have (global) normal crossings if $(D_{\text{red}})_I$ is smooth of codimension $|I|$ in $X$ for all $I = \{i_1, \ldots, i_{|I|}\}$. Notice this is equivalent to requiring the open strata

$$(D_{\text{red}})_I^0 = (D_{\text{red}})_I - \cup_{J \supset I, J \neq I} (D_{\text{red}})_J$$

be smooth of codimension $|I|$.

*Definition* 2.2.   Let $f : X \to S$ be a flat morphism of smooth varieties. Let $Y \subset X$ and $\Sigma \subset S$ be normal crossings divisors, such that $(f^* \Sigma)_{\text{red}} \subset D_{\text{red}}$. The data $\{f : X \to S, Y, \Sigma\}$ are said to be relative normal crossings divisors if for all $I$ there exists a $J$ such that $f(D_{\text{red}})_I \subset (\Sigma_{\text{red}})_J$ and $f : (D_{\text{red}})_I^0 \to (\Sigma_{\text{red}})_J$ is smooth.

LEMMA 2.3.   *Let $\{f : X \to S, Y, \Sigma\}$ be as above. Then $\{f : X \to S, Y, \Sigma\}$ is a relative normal crossings divisor if and only if the sheaf*

$$\Omega^1_{X/S}(\log Y) = \Omega^1_X(\log Y)/f^* \Omega^1_S(\log \Sigma)$$

*is locally free.*

*Proof.* Local freeness is checked in the completion of the local ring of each point. If $\{f : X \to S, Y, \Sigma\}$ is a relative normal crossings divisor, there exist local coordinates $x_j, y_\ell$ on $X$ and $s_i$ on $S$ such that $f$ has local equations $s_i = \Pi_j x_j^{m_{ij}}$, $\Sigma$ has local equation $\Pi_{i \leq r} s_i$, and $Y$ has local equation $\Pi_{i \leq q} s_i \Pi_j x_j \Pi_{\ell \leq p} y_\ell$. Flatness implies that a given $x_j$ appears in at most one $s_i$, so the $ds_i/s_i$ are linearly independent in the fibres of $\Omega^1_S(\log Y)$. A local computation shows that

$$\Omega^1_{X/S}(\log Y) = \Omega^1_X(\log Y)/f^* \Omega^1_S(\log \Sigma)$$

is locally free. The converse is straightforward also.          $\square$

*Remark* 2.4.   The definition of normal crossings divisors does not involve the multiplicities of the components $Y_i$, so these will frequently be ignored. Also, in the relative case, $\Sigma_{\text{red}}$ is determined by $Y$ as it is the image of the union of the components of $Y$ which do not dominate $S$. So given $f : X \to S$ we will simply speak of $Y \subset X$ as a relative normal crossings divisor and use the notation

$$\Omega^1_{X/S}(\log Y) = \Omega^1_X(\log Y)/f^* \Omega^1_S(\log \Sigma).$$

COROLLARY 2.5. *Given $X \xrightarrow{f} S \xrightarrow{g} T$ with divisors $Y \subset X, \Sigma \subset S, \Theta \subset T$ such that both $\{f : X \to S, D, \Sigma\}$ and $\{g : S \to T, \Sigma, \Theta\}$ are relative normal crossings divisors, then $\{g \circ f : X \to T, Y, \Theta\}$ is a relative normal crossings divisor, and the sequence*

$$0 \to f^*\Omega^1_{S/T}(\log \Sigma) \to \Omega^1_{X/T}(\log D) \to \Omega^1_{X/S}(\log D) \to 0$$

*is an exact sequence of locally free sheaves.*

We next recall the theory of relative Chern classes as developed in [27]. Unfortunately, since we need a slight generalization of Saito's results, we are obliged to give some details. Being interested in Chow groups, we will work with $K$-cohomology. We write $\underline{K}_i$ for the $i$th Milnor $K$-sheaf as in Definition 1.6.

PROPOSITION 2.6. *Let $X$ be a smooth variety, and let $Y \subset X$ be a closed subset. Then*

$$H^p_Y(X, \underline{K}_q) \cong H^p_{Y \times \mathbb{A}^n}(X \times \mathbb{A}^n, \underline{K}_q)$$

*for all $p$ and $q$.*

*Proof.* First, by an obvious induction we may suppose $n = 1$. Let $\underline{G}_{q,X}$ denote the Gersten resolution of $\underline{K}_{q,X}$ as in Definitions 1.6, 1.7. Write $\underline{\Gamma}_Y$ for the functor associating to a sheaf its subsheaf of sections with supports in $Y$. If $Y \subset X$ is smooth of codimension $r$,

$$\underline{\Gamma}_Y \underline{G}_{q,X}[r] = \underline{G}_{q-r,Y}; \quad H^p(Y, \underline{K}_q) = H^{p+r}_Y(X, \underline{K}_q).$$

More generally, if $Y$ has pure codimension $r$ with generic points $j : \amalg \mathrm{Spec}\,(F_i) \to Y$ then

$$j_* \oplus_i \underline{G}_{q-r,F_i} \cong \underline{\Gamma}_Y \underline{G}_{q,X} / \varinjlim \underline{\Gamma}_Z \underline{G}_{q,X}$$

where the limit is taken over closed sets $Z \subset Y$ of dimension $< \dim Y$. The proof of the proposition is now by induction on $\dim Y$. If this is 0, then $Y$ is smooth; so by the above it will suffice to show for $F$ a field and $\pi : \mathbb{A}^n_F \to \mathrm{Spec}\,(F)$ the projection, that the natural map

$$\underline{G}_{q,F} \to \pi_* \underline{G}_{q,\mathbb{A}^n_F}$$

is a quasi-isomorphism. This amounts to the assertion

$$H^0(\mathbb{A}^1_F, \underline{K}_q) \cong K_q(F); \quad H^1(\mathbb{A}^1_F, \underline{K}_q) = (0)$$

which in turn follows from the standard exact sequence [1]:

$$(2.1) \quad 0 \longrightarrow K_q(F) \longrightarrow K_q(F(x)) \longrightarrow \oplus_{z \in (\mathbb{A}^1_F)^{(1)}} K_{q-1}(F(z)) \longrightarrow 0.$$

Suppose now we have proved the proposition for $Z \subset X$ of dimension $< \dim Y$. Consider the diagram

(2.2)

$$
\begin{array}{ccccccccc}
0 & \to & \varinjlim \underline{\Gamma}_Z \underline{\underline{G}}_{q,X} & \to & \underline{\Gamma}_Y \underline{\underline{G}}_{q,X} & \to & j_* \oplus_i \underline{\underline{G}}_{q,\mathrm{Spec}\,(F_i)} & \to & 0 \\
& & \downarrow & & \downarrow & & \downarrow & & \\
0 & \to & \varinjlim \pi_* \underline{\Gamma}_{\mathbb{A}^1_Z} \underline{\underline{G}}_{q,\mathbb{A}^1_X} & \to & \pi_* \underline{\Gamma}_{\mathbb{A}^1_Y} \underline{\underline{G}}_{q,\mathbb{A}^1_X} & \to & j_* \pi_* \oplus_i \underline{\underline{G}}_{q,\mathbb{A}^1_{F_i}} & \to & 0.
\end{array}
$$

By induction, the left-hand vertical arrow is a quasi-isomorphism. Also, the complexes $\underline{\underline{G}}_q$ are made up of constant sheaves supported on closed subvarieties, and are therefore acyclic for $j_*$; thus we can think of the right-hand vertical arrow as coming from applying $\mathbb{R}j_*$ to the map (a quasi-isomorphism by the above) $\oplus \underline{\underline{G}}_{q,F_i} \to \pi_* \oplus \underline{\underline{G}}_{q,\mathbb{A}^1_{F_i}}$. Since $\mathbb{R}j_*$ preserves quasi-isomorphisms, it follows that the right-hand map is a quasi-isomorphism, and so the map in the middle is as well. The proposition follows.                                □

A morphism $A \to X$ is called an affine bundle of dimension $n$ if Zariski-locally on $X$, $A \cong \mathbb{A}^n_X$. (We do not require the transition maps to be linear.) The following is proved just as above:

COROLLARY 2.7.  *For $X$ smooth, and $A \to X$ an affine bundle,*

$$
H^p(X, \underline{\underline{K}}_q) \cong H^p(A, \underline{\underline{K}}_q)
$$

*for all $p$ and $q$.*

Suppose now that $X$ is smooth as above, and $Y = \cup_{i \in \mathcal{A}} Y_i$ is a normal crossings divisor in $X$. Let $\pi : V \to X$ be a vector bundle, and let $\Delta \subset \pi^{-1}(Y)$ be a subscheme. We assume $\Delta = \cup \Delta_i$ with $\Delta_i \subset \pi^{-1}(Y_i)$ and denote by $Y_I, \Delta_I$ for $I \subset \mathcal{A}$ the intersections $\cap_{i \in I} Y_i, \cap_{i \in I} \Delta_i$. If we order the index set $\mathcal{A}$, we can define as above two sorts of relative $K$-cohomology:

$$
\begin{aligned}
(2.3) \quad \mathbb{H}^*(X, \underline{\underline{K}}_{q,X,Y}) &:= \mathbb{H}^*(X, \underline{\underline{K}}_{q,X} \to \underline{\underline{K}}_{q,Y}), \\
\mathbb{H}^*(X, \underline{\underline{K}}_{q,X,Y_\bullet}) &:= \mathbb{H}^*(X, \underline{\underline{K}}_{q,X} \to \oplus_i \underline{\underline{K}}_{q,Y_i} \to \oplus_{i<j} \underline{\underline{K}}_{q,Y_{i,j}} \to \dots)
\end{aligned}
$$

and a map between them

$$
(2.4) \qquad \mathbb{H}^*(X, \underline{\underline{K}}_{q,X,Y}) \to \mathbb{H}^*(X, \underline{\underline{K}}_{q,X,Y_\bullet}).
$$

Similarly, we can define

$$
(2.5) \qquad \mathbb{H}^*(V, \underline{\underline{K}}_{q,V,\Delta}) \to \mathbb{H}^*(V, \underline{\underline{K}}_{q,V,\Delta_\bullet}).
$$

PROPOSITION 2.8. *Let notation be as above, and assume that $\Delta_I = \cap_{i \in I} \Delta_i$ is a nonempty affine bundle over $Y_I$ for all $I \subset \mathcal{A}$ with $|I| \leq p$. Then there is an isomorphism*

$$\pi^* : \mathbb{H}^m(X, \underline{K}_{q,X,Y_\bullet}) \tilde{\rightarrow} \mathbb{H}^m(V, \underline{K}_{q,V,\Delta_\bullet})$$

*for all $m < p$.*

*Proof.* There is a spectral sequence

$$E_1^{a,b}(X, Y_\bullet) = \oplus_{\#I=a} H^b(Y_I, \underline{K}_q) \Rightarrow \mathbb{H}^*(X, \underline{K}_{q,X,Y_\bullet}),$$

and by Corollary 2.7 we have $\pi^* : E_1^{a,b}(X, Y_\bullet) \cong E_1^{a,b}(V, \Delta_\bullet)$ whenever $a \leq p$. $\qquad\square$

Now let $W$ be an algebraic cycle of codimension $r$ on $V$, and assume the support $|W|$ does not meet $\Delta$. The local $K$-cohomology carries a cycle class, so that

(2.6)
$$[W] \in H^r_{|W|}(V, \underline{K}_r) \cong H^r_{|W|}(V, \underline{K}_{r,V,\Delta_\bullet}) \rightarrow H^r(V, \underline{K}_{r,V,\Delta_\bullet}) \cong \mathbb{H}^r(X, \underline{K}_{r,X,Y_\bullet}).$$

In the next examples, we define the class appearing on the right-hand side of the Riemann-Roch formula 0.1.

*Example* 2.9 (T. Saito [27]). Suppose there exist vector bundle surjections

$$\phi_i : V|_{Y_i} \twoheadrightarrow \mathcal{O}_{Y_i}$$

which are independent in the sense that for any $I$, the map

$$\oplus_{i \in I} \phi_i : V|_{Y_I} \twoheadrightarrow \oplus_{i \in I} \mathcal{O}_{Y_I}$$

is surjective. Define $\Delta_i := \phi_i^{-1}(1) \subset V|_{Y_i}$. Take $W = 0$-section $\subset V$. T. Saito defines, for $d = \text{rk}V$,

$$c_d(V, \phi) := [W] \in \mathbb{H}^d(X, \underline{K}_{d,X,Y_\bullet}).$$

Note that, in Saito's case, the $\Delta_i$ meet properly.

*Example* 2.10. Suppose there exists a single surjective map

$$\phi : V|_Y \twoheadrightarrow \mathcal{O}_Y.$$

Define $\phi_i = \phi|_{Y_i}$ and $\Delta_i = \phi_i^{-1}(1)$. In this case the $\Delta_I \rightarrow Y_I$ are all affine bundles of fibre dimension $d - 1$, but we may still define

$$c_d(V, \phi) \in \mathbb{H}^d(X, \underline{K}_{d,X,Y_\bullet}).$$

Next we check the compatibility of the relative classes with exact sequences, which we will need later in various reduction steps.

PROPOSITION 2.11.    *Let $X$ be smooth and $Y = \cup_{i \in \mathcal{A}} Y_i$ a reduced normal crossings divisor as above. Let*

$$0 \longrightarrow V' \longrightarrow V \longrightarrow V'' \longrightarrow 0$$

*be an exact sequence of vector bundles on $X$ of ranks $d', d, d''$ respectively. Let $\phi_i : V|_{Y_i} \twoheadrightarrow \mathcal{O}_{Y_i}$ be a partial trivialization, and assume $\Delta_I \to Y_I$ is an affine bundle for all $I \subset \mathcal{A}$. Suppose there is a decomposition $\mathcal{A} = \mathcal{A}' \amalg \mathcal{A}''$ such that for $i \in \mathcal{A}''$, $\phi_i(V') = (0)$; now there is an induced partial trivialization of $V''$ over $Y'' = \cup_{i \in \mathcal{A}''} Y_i$. Assume also, for $I \subset \mathcal{A}'$, that $V' \cap \Delta_I \to Y_I$ is an affine bundle. Then the relative top Chern classes*

$$c_{d'}(V', \{\phi_i|_{V'}\}_{i \in \mathcal{A}'}), \ c_d(V, \{\phi_i\}_{i \in \mathcal{A}}), \ c_{d''}(V'', \{\phi_i\}_{i \in \mathcal{A}''})$$

*are defined (in the groups $\mathbb{H}^*(X, \underline{\underline{K}}_{*,X,Z_\bullet})$ with $Z = Y', Y, Y''$ respectively), and*

$$c_d(V, \{\phi_i\}_{i \in \mathcal{A}}) = c_{d'}(V', \{\phi_i|_{V'}\}_{i \in \mathcal{A}'}) c_{d''}(V'', \{\phi_i\}_{i \in \mathcal{A}''}).$$

*Proof.* We leave for the reader the construction of a product

$$(2.7) \qquad \mathbb{H}^a(X, \underline{\underline{K}}_{c,X,Y'_\bullet}) \otimes \mathbb{H}^b(X, \underline{\underline{K}}_{d,X,Y''_\bullet}) \to \mathbb{H}^{a+b}(X, \underline{\underline{K}}_{c+d,X,Y_\bullet})$$

compatible with augmentation to the usual (nonrelative) $K$-cohomology. This can be done, for example, if we use a variant of the usual cochain product

$$(xy)(i_0, \ldots, i_{a+b}) = x(i_0, \ldots, i_a) \cdot y(i_a, \ldots, i_{a+b}).$$

Having done this, suppose $W', W''$ are cycles on $X$ of codimensions $a, b$ disjoint from $Y', Y''$ and meeting properly. Let $W = W' \cdot W''$. One has a product on local cohomology

$$H^a_{|W'|}(X, \underline{\underline{K}}_c) \times H^b_{|W''|}(X, \underline{\underline{K}}_d) \to H^{a+b}_{|W|}(X, \underline{\underline{K}}_{c+d})$$

which is compatible with the cycle classes. Since these local cohomology groups are isomorphic to the corresponding local relative cohomology groups, one gets in this case that the product on relative cohomology is compatible with the cycle product. Another advantage of the cycle class construction in local cohomology is that it can be done locally on $X$ and the classes glued.

We will apply the above discussion with $X$ replaced by the total space of the vector bundle $V$. Replacing $X$ by an affine torseur over $X$ (which does not affect the relative $K$-cohomology because all the $Y_I$ are smooth) we may assume $V = V' \oplus V''$. Let $p' : V \to V'$ and $p'' : V \to V''$ be the projection maps. We define

$$W' = p'^*(0\text{-section of } V'); \quad W'' = p''^*(0\text{-section of } V'')$$

so that $W = W' \cdot W''$ is the zero section of $V$. Note finally that the composition

$$\mathbb{H}^a(X, \underline{\underline{K}}_{c,X,Y_\bullet}) \cong \mathbb{H}^a(V', \underline{\underline{K}}_{c,V',\Delta_\bullet}) \overset{p^*}{\underset{\cong}{\cong}} \mathbb{H}^a(V, \underline{\underline{K}}_{c,V,p'^{-1}(\Delta_\bullet)})$$
$$\cong \mathbb{H}^a(X, \underline{\underline{K}}_{c,X,Y'_\bullet})$$

is the identity. The rest of the argument is straightforward, given the identification of the top Chern class with the class of the 0-section. $\square$

*Example* 2.12. Let $\{f : X \to S, Y, \Sigma\}$ be a relative normal crossings divisor as in Definition 2.2, and consider the exact sequence

$$0 \to f^*\Omega^1_S(\log \Sigma) \to \Omega^1_X(\log Y) \to \Omega^1_{X/S}(\log Y) \to 0.$$

Write $Y = \cup_{i \in \mathcal{A}} Y_i$ and let $\mathcal{A}' = \{i \in \mathcal{A} \,|\, D_i \subset f^{-1}\Sigma\}$, $\mathcal{A}'' = \mathcal{A} - \mathcal{A}'$, $Y' = \cup_{i \in \mathcal{A}''} D_i$, $Y'' = \cup_{i, \mathcal{A}''} Y_i$. Define a partial trivialization

$$\mathrm{res}_i = \mathrm{res}_{Y_i} : \Omega^1_X(\log Y)|_{Y_i} \twoheadrightarrow \mathcal{O}_{Y_i}.$$

If $f^*(\Sigma)$ is reduced we are in the situation of Proposition 2.11, and we may conclude

$$c_{\dim S}(f^*\Omega^1_S(\log \Sigma), \mathrm{res}'_Y) \cdot c_{\dim(X/S)}(\Omega^1_{X/S}(\log Y), \mathrm{res}''_Y)$$
$$= c_{\dim X}(\Omega^1_X(\log Y), \mathrm{res}_Y).$$

Note however, that if $f^*(\Sigma)$ is not reduced, the induced partial trivialization of the left-hand bundle is not the pullback of the partial trivialization on $S$. In this case it can happen that for some $I \subset \mathcal{A}'$ we have $\Delta_I \cap f^*\Omega^1_S(\log \Sigma) = \emptyset$. Define a modified partial trivialization

$$\rho_i = \begin{cases} \mathrm{res}_i & i \in \mathcal{A}'' \\ \mathrm{ord}_{Y_i}(f^*(\Sigma))^{-1} \cdot \mathrm{res}_i & i \in \mathcal{A}'. \end{cases}$$

This partial trivialization is compatible with the pullback of res on $\Omega^1_S(\log \Sigma)$. Omitting a straightforward verification of contravariant functoriality, we deduce for a suitable error term $\varepsilon$

(2.8) $$f^*(c_{\dim S}(\Omega^1_S(\log \Sigma), \mathrm{res}_\Sigma) \cdot c_{\dim X/S}(\Omega^1_{X/S}(\log Y), \mathrm{res}_{D''})$$
$$= c_{\dim X}(\Omega^1_X(\log Y), \rho) = c_{\dim X}(\Omega^1_X(\log Y), \mathrm{res}_Y) + \varepsilon.$$

The error term $\varepsilon$ in the above example has been calculated by Saito ([27, Prop. 1]). It will turn out to be inoffensive for our purposes. To see this we need to look more closely at the product on the right-hand side of Riemann-Roch 0.15. As in Section 1, we work with the complexes

(2.9) $$\Omega\underline{\underline{K}}_{n,X}(\log D) := \tilde{j}_*\underline{\underline{K}}_{n,X-D} \xrightarrow{d\log} \Omega^n_X(\log D)$$
$$\to \Omega^{n+1}_X(\log D) \to \ldots \to \Omega^{2n-1}_X(\log D).$$

Here $\tilde{j} : X - D \to X$ is the inclusion. We write $\Omega\underline{\underline{K}}_{n,X}$ when $D = \emptyset$.

We consider a normal crossings divisor $Y = W + Z$. In the application, $\{f : X \to S, Y, \Sigma\}$ will be a relative normal crossings divisor, with $W = f^{-1}(\Sigma)_{\mathrm{red}}$. We define a pairing of complexes in the derived category

$$(2.10) \qquad \Omega\underline{\underline{K}}_{n,X}(\log Y) \times \underline{\underline{K}}_{d,X,Z_\bullet} \to \Omega\underline{\underline{K}}_{d+n,X}(\log W).$$

Let $\tilde{D}_I := D_I - \cup_K D_K$ where $K \supset I, K \neq I$. Here $D$ denotes either $Y$, or $W$ or $Z$. Let

$$\tilde{j}_p : \tilde{D}^{(p)} \quad := \quad \amalg_{\#I=p}\tilde{D}_I \to X,$$
$$j_p : D^{(p)} \quad := \quad \text{normalization of } \amalg_{\#I=p} D_I \to X$$

be the inclusion. The following double complex $C^{a,b}, a, b \geq 0$, is quasi-isomorphic to $\Omega\underline{\underline{K}}_{d+n,X}(\log W)$:

$$(2.11)$$

Indeed, we will show in the following lemma that the column starting with $\tilde{j}_*\underline{\underline{K}}_{d+n,X-Y}$ is quasi-isomorphic to $\tilde{j}_*\underline{\underline{K}}_{d+n,X-W}$. The standard residue sequence shows that the column starting with $\Omega^i_X(\log Y)$ is quasi-isomorphic to $\Omega^i_X(\log W)$.

LEMMA 2.13. *Let $X = \mathrm{Spec}(R)$ be the spectrum of a local ring on a smooth variety. Let $D = \bigcup_{i=1}^r D_i \subset X$ be a normal crossings divisor, and let $U = X - D, U_s = X - \bigcup_{i=s+1}^r D_i$. Then the Gersten complex $G_q(U) := H^0(U, \underline{G}_{q,X})$ is a resolution of $K_q(U) := \Gamma(U, \underline{K}_q)$ (cf. Definitions 1.6, 1.7). When $D_I = \bigcap_{i \in I} D_i$ and $\tilde{D}_I = D_I - \bigcup_{J \subsetneq I} D_J$, there is an exact sequence*

$$0 \to K_q(U_s) \to K_q(U) \to \bigoplus_{i=1}^s K_{q-1}(\tilde{D}_i) \to \bigoplus_{1 \leq i < j \leq s} K_{q-2}(\tilde{D}_{\{i,j\}}) \to \dots .$$

*Proof.* We will prove the statement about the Gersten complex by induction on $r$. If $r = 0$ this is proved in [25]. Assume $r \geq 1$ and the lemma holds for $r - 1$. Write

$$T = D_1; \ E = T \cap (\cup_{i=2}^r D_i); \ V = T - E = \tilde{D}_1; \ U' = U_1.$$

Let $F = k(X)$ and $L = k(T)$ be the function fields. Consider the diagram

$$
\begin{array}{ccccccccc}
& & & & 0 & & 0 & & \\
& & & & \downarrow & & \downarrow & & \\
& & & & K_{q-1}(k(T)) & \xrightarrow{b} & \coprod_{v \in V^1} K_{q-2}(k(v)) & & \\
& & 0 & & \downarrow & & \downarrow & & \\
& & \downarrow & & & & & & \\
0 & \to & K_q(U') & \to & K_q(k(X)) & \to & \coprod_{U'^1} K_{q-1}(k(x)) & \to & \coprod_{U'^2} K_{q-2}(k(x)) \\
& & \downarrow & & \| & & \downarrow & & \downarrow \\
0 & \to & K_q(U) & \xrightarrow{c} & K_q(k(X)) & \xrightarrow{d} & \coprod_{U^1} K_{q-1}(k(x)) & \to & \coprod_{U^2} K_{q-2}(k(x)) \\
& & {\scriptstyle a}\downarrow & & & & \downarrow & & \downarrow \\
& & K_{q-1}(V) & = & \ker(b) & & 0 & & 0
\end{array}
$$

The middle row is exact by induction, and the top row is exact except that $\ker(b) \cong K_{q-1}(V)$. Also the map $c$ is injective, and $\text{image}(c) = \ker(d)$. The columns except possibly the first are also exact. One then proves the surjectivity of $a$ and the exactness of the first column by a diagram chase. Another chase gives exactness of the resolution of $K_q(U)$.

The second part of the lemma is now proved by induction on $s$. For $s = 1$ it is the right-hand column of the above diagram. Assume the assertion for $s - 1$, and consider the diagram

$$
\begin{array}{ccccccccc}
& & & & 0 & & 0 & & \\
& & & & \downarrow & & \downarrow & & \\
& & & & K_{q-1}(\tilde{D}_s)) & \xrightarrow{b} & \bigoplus_{1 \le i < s} K_{q-2}(\tilde{D}_{\{i,s\}}) \to \cdots & & \\
& & 0 & & \downarrow & & \downarrow & & \\
& & \downarrow & & & & & & \\
0 & \to & K_q(U_s) & \to & K_q(U) & \to & \bigoplus_{i=1}^{i=s} K_{q-1}(\tilde{D}_i) & \to & \bigoplus_{1 \le i < j \le s} K_{q-2}(\tilde{D}_{\{i,j\}}) \to \cdots \\
& & \downarrow & & \| & & \downarrow & & \downarrow \\
0 & \to & K_q(U_{s-1}) & \xrightarrow{c} & K_q(U) & \xrightarrow{d} & \oplus_{i=1}^{s-1} K_{q-1}(\tilde{D}_i) & \to & \oplus_{1 \le i < j \le s-1} K_{q-2}(\tilde{D}_{\{i,j\}}) \to \cdots \\
& & {\scriptstyle a}\downarrow & & & & \downarrow & & \downarrow \\
& & \ker(b) & & & & 0 & & 0
\end{array}
$$

By induction, $\ker(b) \cong K_{q-1}(D_s - \bigcup_{t > s} D_t)$ and (from the previous diagram) the map $a$ is surjective. The bottom and top rows are exact (except $\ker(b) \ne (0)$). Again a diagram chase shows the middle row is exact. $\square$

Returning to the construction of the product, we order the index set $\mathcal{A}$ of components of $Y$. Thus, for

$$I = \{i_1 < \ldots < i_r\} \subset \mathcal{A}$$

we can define an iterated tame symbol and an iterated residue

$$(2.12) \quad t_I \quad = \quad t_{i_r} \circ t_{i_{r-1}} \circ \ldots \circ t_{i_1} : \tilde{j}_* \underline{\underline{K}}_{p,X-Y} \to \tilde{j}_{p-r,*} \underline{\underline{K}}_{p-r,Z^{(r)}}$$

$$\mathrm{res}_I \quad = \quad \mathrm{res}_{i_r} \circ \mathrm{res}_{i_{r-1}} \circ \ldots \circ \mathrm{res}_{i_1} : \Omega^s_X(\log Y) \to \Omega^{s-r}_{Z^{(r)}}(\log Z^{(r+1)}).$$

We define pairings ($r = \#I, s \geq 1$)

$$(2.13) \quad \tilde{j}_* \underline{\underline{K}}_{n,X-Y} \times \underline{\underline{K}}_{l,X,Z_\bullet} \to C^{0,\cdot}; \quad a \times b_I \mapsto t_I(a) \cdot b_I \in \tilde{j}_{r,*} \underline{\underline{K}}_{n-r,Z^{(r)}},$$

$$\Omega^{n+s-1}_X(\log Y) \times \underline{\underline{K}}_{q,X,Z_\bullet} \to C^{s,\cdot}; \quad a \times b_I \mapsto \mathrm{res}_I(a) \wedge d\log(b_I).$$

These induce the desired pairing (2.10).

*Remark* 2.14. Given $(E, \nabla) \in \mathcal{G}(X, \log Y)$, its Chern-Simons classes $\mathrm{N}w_n((E, \nabla))$ lie in

$$H^{2n}_{\mathrm{CS}}(X(\log Y)) = \Gamma(X, \underline{H}^n(\Omega \underline{K}_{n,X}(\log Y)))$$

(see Corollary 1.10). On the other hand, if $\{f : X \to S, Y, \Sigma\}$ is a relative normal crossings divisor (see Definition 2.2) of relative dimension $d$, then there is the relative top Chern class

$$c_d(\Omega^1_{X/S}(\log Y), \mathrm{res}_Z) \in \mathbb{H}^d(X, \underline{\underline{K}}_{d,X,Z_\bullet})$$

where $Y = f^{-1}(\Sigma)_{\mathrm{red}} + Z$ (see 2.9). The pairing (2.10) induces, for $n \geq 1$,

$$\underline{H}^n(\Omega \underline{K}_{n,X}(\log Y)) \times \mathbb{H}^d(X, \underline{\underline{K}}_{d,X,Z_\bullet}) \overset{\cdot}{\to} \underline{\underline{H}}^{d+n}(\Omega \underline{\underline{K}}_{d+n,X}(\log f^{-1}(\Sigma)_{\mathrm{red}})),$$

defining

$$(2.14) \quad c_d(\Omega^1_{X/S}(\log Y), \mathrm{res}_Z) \cdot \mathrm{N}w_n((E, \nabla)) \in H^{2(d+n)}_{\mathrm{CS}}(X(\log(f^{-1}\Sigma)_{\mathrm{red}}))$$

$$= \Gamma(X, \underline{\underline{H}}^{d+n}(\Omega \underline{\underline{K}}_{d+n,X}(\log f^{-1}(\Sigma)_{\mathrm{red}}))).$$

We can now show that the error term $\varepsilon$ from example 2.12 does not affect our Riemann-Roch calculations. We consider a relative normal crossings situation $\{f : X \to S, Y, \Sigma\}$ with $d = \dim X/S$. Saito's calculations imply for $n_i \geq 1$

$$(2.15) \quad c_d(\Omega^1_{X/S}(\log Y), \mathrm{res}_i) - c_d(\Omega^1_{X/S}(\log Y), n_i \cdot \mathrm{res}_i)$$

$$\in \mathrm{image}\Big(K_1(\mathbb{Q}) \otimes H^d(X, \underline{\underline{K}}_{d-1,X,Z_\bullet}) \to H^d(X, \underline{\underline{K}}_{d,X,Z_\bullet})\Big).$$

LEMMA 2.15. *The pairing*

$$K_1(\mathbb{Q}) \otimes H^d(X, \underline{\underline{K}}_{d-1,X,Z_\bullet}) \otimes \underline{H}^n(\Omega \underline{K}_{n,X}(\log Y))$$

$$\to \underline{\underline{H}}^{d+n}(\Omega \underline{\underline{K}}_{d+n,X}(\log f^{-1}(\Sigma)_{\mathrm{red}}))$$

*is zero.*

*Proof.* This pairing can be factored

$$K_1(\mathbb{Q}) \otimes H^d(X, \underline{\underline{K}}_{d-1,X,Z_\bullet}) \otimes \underline{\underline{H}}^n(\Omega\underline{\underline{K}}_{n,X}(\log Y))$$

$$\rightarrow K_1(\mathbb{Q}) \otimes \underline{\underline{H}}^{d+n}(\Omega\underline{K}_{d+n-1,X}(\log f^{-1}(\Sigma)_{\text{red}}))$$

$$\rightarrow \underline{\underline{H}}^{d+n}(\Omega\underline{K}_{d+n,X}(\log f^{-1}(\Sigma)_{\text{red}})).$$

For $a \in \mathbb{Q}^\times$, the second arrow comes from the map on complexes

$$a \cdot ? : \Omega\underline{K}_{d+n-1,X}(\log f^{-1}(\Sigma)_{\text{red}}) \rightarrow \Omega\underline{K}_{d+n,X}(\log f^{-1}(\Sigma)_{\text{red}})$$

which is multiplication by $a$ on the $K$-sheaf and zero on the differentials. This induces 0 on $\underline{\underline{H}}^{d+n}$ because $n \geq 1$. $\square$

As a consequence we have:

PROPOSITION 2.16. *Let $\{f : X \rightarrow S, Y, \Sigma\}$ be a relative normal crossings divisor with $\dim X/S = d$, with $Y = f^{-1}(\Sigma_{\text{red}}) + Z$. Let $(E, \nabla) \in \mathcal{G}(X, \log Y)$. Let $n_i \geq 1$ be a collection of multiplicities. Then*

$$(2.16) \quad c_d(\Omega^1_{X/S}(\log Y), \text{res}_i) \cdot \text{Nw}_n((E, \nabla))$$

$$= c_d(\Omega^1_{X/S}(\log Y), n_i \cdot \text{res}_i) \cdot \text{Nw}_n((E, \nabla))$$

$$\in H^{2(d+n)}_{\text{CS}}(X(\log(f^{-1}(\Sigma)_{\text{red}})).$$

Finally, we define a transfer map

$$(2.17) \quad f_* : H^{2(d+p)}_{\text{CS}}(X(\log(f^{-1}(\Sigma)_{\text{red}}))) \rightarrow H^{2p}_{\text{CS}}(S(\log \Sigma))$$

using Cousin complexes. Recall [18], for $\mathcal{F}$ an abelian sheaf on $X$, the complex Cousin($\mathcal{F}$) is given by

$$\coprod_{x \in X^0} i_{x*}\mathcal{F}_x \rightarrow \coprod_{x \in X^1} i_{x*}H^1_{\{x\}}(\mathcal{F}) \rightarrow \coprod_{x \in X^2} i_{x*}H^2_{\{x\}}(\mathcal{F}) \rightarrow \dots .$$

Here $X^r$ is the set of points of codimension $r$ in the scheme $X$, $H^i_{\{x\}}(\mathcal{F})$ denotes the $i^{\text{th}}$ local cohomology of $\mathcal{F}$ with supports in $\{x\}$, and $i_{x*}A$ for an abelian group $A$ is the direct image on $X$ of the constant sheaf $A_{\overline{\{x\}}}$ on the Zariski closure $\overline{\{x\}}$ of the point $x$. For $X$ smooth over a field, the Cousin complex is a resolution of $\mathcal{F}$ for $\mathcal{F} = \underline{K}_n$ with the Definition 1.6, (in which case, $H^i_{\{x\}}(\mathcal{F}) \cong K_{n-i}(k(x))$), and also for $\mathcal{F}$ coherent and locally free.

PROPOSITION 2.17.    *Let $f : X \to S$ be a proper morphism of smooth varieties. Let $\Sigma \subset S$ be a normal crossings divisor, and assume $W = f^{-1}(\Sigma)_{\mathrm{red}}$ is a normal crossings divisor as well. Let $j : X - W \to X$, $i : S - \Sigma \to S$ be the open embeddings and $d = \dim X - \dim S$.   Then there exist transfer maps*

$$(2.18) \qquad \mathrm{tr} : f_* \mathrm{Cousin}(\Omega_X^r(\log W)) \quad \to \quad \mathrm{Cousin}(\Omega_S^{r-d}(\log \Sigma))[-d]$$

$$\mathrm{tr} : f_* \mathrm{Cousin}(j_* \underline{\underline{K}}_{n,X-W}) \quad \to \quad \mathrm{Cousin}(i_* \underline{\underline{K}}_{n-d,S-\Sigma})[-d].$$

*These maps are compatible with $d$ and $d \log$.*

*Proof.*  For the $K$-sheaves, the Cousin complex coincides with the Gersten resolution. Further, using Lemma 2.13 above, we get that $j_*(\mathrm{Cousin}(\underline{\underline{K}}_{n,X-W}))$ is a resolution of $j_* \underline{\underline{K}}_{n,X-W}$.

We define the transfer on the double complex

$$(2.19) \qquad j_* \mathrm{Cousin}(\underline{\underline{K}}_{n,X-W}) \xrightarrow{d \log} \mathrm{Cousin}(\Omega_X^n(\log W))$$
$$\to \ldots \to \mathrm{Cousin}(\Omega_X^{2n-1}(\log W)).$$

Taking $x \in X^r$ in $j_* \mathrm{Cousin} \underline{\underline{K}}_{n,X-W}$, the transfer map

$$H_{\{x\}}^r(\underline{\underline{K}}_{n,X-W}) \cong K_{n-r}(k(x)) \to K_{n-r}(k(f(x))) \cong H_{\{f(x)\}}^{r-d}(\underline{\underline{K}}_{n-d,S-\Sigma})$$

is the trace if $[k(x) : k(f(x))] < \infty$ and is zero otherwise. For details on this $K$-theoretic trace, see [4], [20].

The construction in the case of differential forms is built around an iterated residue. When $r = 0$ it is just the trace on differential forms from the function field on $X$ to the function field on $S$. This trace carries forms with log poles on $f^{-1}(\Sigma)$ to forms with log poles on $\Sigma$. It is zero if $[k(X) : k(S)] = \infty$. Suppose next that $r = 1$ and locally near $x$ the subscheme $\overline{\{x\}} = T$ is defined by $t = 0$ on $X$, and $x \notin W$. If $[k(x) : k(f(x))] = \infty$ the transfer $H_{\{x\}}^1(\Omega_X^n) \to H_{\{f(x)\}}^{1-d}(\Omega_S^{n-d})$ is zero. Assume $[k(x) : k(f(x))] < \infty$, so that $d \leq 1$. If $d = 0$ then $f(x)$ is a codimension 1 point on $S$. Let $s = 0$ be a local defining equation. The transfer is defined to be the composition

$$H_{\{x\}}^1(\Omega_X^n) \quad = \quad \Omega_{X,x}^n[t^{-1}] \Big/ \Omega_{X,x}^n \hookrightarrow \Omega_{X,x}^n[f^*(s)^{-1}] \Big/ \Omega_{X,x}^n$$
$$\xrightarrow{\mathrm{Tr}} \Omega_{S,f(x)}^n[s^{-1}] \Big/ \Omega_{S,f(x)}^n.$$

If $d = 1$, then $f(x)$ is a codimension 0 point on $S$. One has

$$H_{\{x\}}^1(\Omega_X^n) \quad \twoheadrightarrow \quad \Omega_{X,x}^n[t^{-1}] \Big/ (\Omega_{X,x}^n + f^* \Omega_{S,f(x)}^n)$$
$$\cong \Omega_{X,x}^{n-1} \otimes \left( \Omega_{X/S,x}^1[t^{-1}] \Big/ \Omega_{X/S,x}^1 \right) \xrightarrow{\mathrm{res}_x} \Omega_{k(x)}^{n-1} \xrightarrow{\mathrm{tr}} \Omega_{k(f(x))}^{n-1},$$

which is the transfer in this case.

Suppose now $r > 1$ and $x \notin W$. Finiteness of $[k(x) : k(f(x))]$ implies $f(x) : s_1 = \ldots = s_{r-d} = 0$. Let $t_i = f^* s_i$ for $i \le r-d$, and choose $t_{r-d+1}, \ldots, t_r$ such that $t_1 = \ldots = t_r = 0$ locally defines some multiple of $x$. The desired transfer map

$$
\begin{aligned}
H^r_{\{x\}}(\Omega^n_X) \quad &\cong \quad \Omega^n_{X,x}[(t_1 \cdots t_r)^{-1}] \Big/ \bigg( \sum_i \Omega^n_{X,x}[(t_1 \cdots \hat{t}_i \cdots t_r)^{-1}] \bigg) \\
&\to \Omega^{n-d}_{S,f(x)}[(s_1 \cdots s_{r-d})^{-1}] \Big/ \bigg( \sum_i \Omega^{n-d}_{S,f(x)}[(s_1 \cdots \hat{s}_i \cdots s_{r-d})^{-1}] \bigg) \\
&\cong \quad H^{r-d}_{\{f(x)\}}(\Omega^{n-d}_S)
\end{aligned}
$$

is defined by an iterated residue $\Omega^d_{X/S,x}[(t_{r-d+1} \cdots t_r)^{-1}] \to k(x)$. Details are omitted.

Finally, suppose $x$ is a codimension $r$ point on $X$ which lies on $W = f^{-1}(\Sigma)$. Write $W^{(i)}$ for the normalized $i$-fold intersection of components, and let $\{x_i\} \subset W^{(i)}$ be the set of points lying over $x$. Similarly, suppose $f(x)$ lies on $\Sigma$. We may calculate the local cohomology of the log forms using $E_1$ spectral sequences associated to the weight filtrations. One gets a diagram

$$
\begin{array}{ccccccccc}
\cdots \longrightarrow & H^{r-1}_{\{x_1\}}(\Omega^{n-1}_{W^{(1)}}) & \longrightarrow & H^r_{\{x\}}(\Omega^n_X) & \longrightarrow & H^r_{\{x\}}(\Omega^n_X(\log W)) & \longrightarrow & 0 \\
& \downarrow a & & \downarrow b & & \downarrow c & & \\
\cdots \longrightarrow & H^{r-1-d}_{\{f(x_1)\}}(\Omega^{n-1}_{\Sigma^{(1)}}) & \longrightarrow & H^{r-d}_{\{f(x)\}}(\Omega^n_S) & \longrightarrow & H^{r-d}_{\{f(x)\}}(\Omega^n_S(\log \Sigma)) & \longrightarrow & 0.
\end{array}
$$

The trace maps labelled $a$ and $b$ are constructed as above, and the rows are exact by purity, and so the desired trace map $c$ is defined. (The point here is that $H^i_{\{x\}}(X, \mathrm{gr}^p_W(\Omega^n(\log W)) = (0)$ unless $i = r - p$.) $\qquad \square$

*Example* 2.18. The terms of the Cousin complex are constant sheaves supported on closed subsets and hence are acyclic for the Zariski topology. We obtain (with notation as above) a transfer map

$$
\mathrm{tr}_{X/S} : \mathbb{R}f_* \Omega \underline{\underline{K}}_{n,X}(\log W) \to \Omega \underline{\underline{K}}_{n-d,S}(\log \Sigma)[-d].
$$

In particular, we obtain a map

$$
f_* : H^{2n}_{\mathrm{CS}}(X(\log f^{-1}(\Sigma))) \to H^{2(n-d)}_{\mathrm{CS}}(S(\log \Sigma)).
$$

In combination with the product (2.10), we have now defined the right-hand side of our Riemann-Roch theorem

$$
f_* \bigg( c_d(\Omega^d_{X/S}(\log Y), \mathrm{res}_Z) \cdot \mathrm{Nw}_n(E, \nabla) \bigg).
$$

## 3. Reductions

Our objective in this section is to reduce the Riemann-Roch theorem to the case of a bundle with log poles on $\mathbb{P}^1_F$ where $F = k(S)$ is a function field. We begin with a relative normal crossings divisor $\{f : X \to S, Y, \Sigma\}$ as in Definition 2.2 with $\dim X/S = d$. Let $(E, \nabla)$ be a connection with logarithmic poles along $Y$ (see Definition 1.2).

*Definition* 3.1. The curvature form of $\nabla$ is *basic*, if

$$(3.1) \quad \nabla^2 \in \operatorname{Hom}_{\mathcal{O}_X}(E, f^*\Omega^2_S(\log \Sigma) \otimes E) \subset \operatorname{Hom}_{\mathcal{O}_X}(E, \Omega^2_X(\log Y) \otimes E).$$

Of course, this is satisfied when $\nabla$ is integrable. It also holds for tensor connections

$$(3.2) \qquad (E, \nabla) = f^*(M, \nabla_M) \otimes (N, \nabla_N), \ (\nabla_N)^2 = 0$$

where $(M, \nabla_M)$ is a connection on $S$ with logarithmic poles along $\Sigma$ and $(N, \nabla_N)$ is an integrable connection on $X$ with logarithmic poles along $Y$. If $M$ is locally free, then the projection formula implies that the Riemann-Roch formula (0.15) for $\nabla$ is a formal consequence of the Riemann-Roch formula for $\nabla_N$. As an example, one can consider $f : X = Z \times S \to S, Y = \Sigma = \emptyset$. Then $(E, \nabla)$ admits a filtration by subbundles with connection $(E_i, \nabla|E_i)$, such that the graded pieces $(E_i/E_{i-1}, \nabla)$ are tensor connections with $\nabla_N$ a flat connection coming from $Z$ ([11]). Thus in this case, the Riemann-Roch formula (0.15) is trivial, even with coefficients in AD, and the result of [13] is of no interest.

On the other hand, poles introduce some flexibility and it is interesting to consider the basic curvature form condition on the connections $(\oplus^N_1 \mathcal{O}_X, \nabla)$ treated in Section 4, where $X = \mathbb{P}^1_S$, $S = \operatorname{Spec} F$, $F$ is a function field, and $\nabla = \Phi + \sum A_i d \log(z - a_i)$. Here $\Phi$ is an $N \times N$-matrix of one-forms on $S$ relative to $k$, $A_i$ is an $N \times N$-matrix with coefficients in $F$, and $a_i : S \to \mathbb{P}^1$ is a section. Then the condition is equivalent to the system of 1-st order differential equations

$$(3.3) \qquad dA_i = [\Phi, A_i] - \sum_{j \neq i}[A_i, A_j] d \log(a_i - a_j).$$

For example, if $A_i = \lambda_i \operatorname{Id}_{N \times N}$, where $\lambda_i \in k$, then the basic curvature form condition is fulfilled. The curvature of such a connection satisfies $\nabla^2 = f^*(d\Phi - \Phi \wedge \Phi)$, but the Chern-Simons classes are not in general pulled back from the base.

Condition 3.1 implies that the relative connection

$$\nabla_{X/S} : E \to \Omega^1_{X/S}(\log Y) \otimes E$$

is integrable, and thus one has coherent Gauss-Manin sheaves:

$$R^i f_*(\Omega^*_{X/S}(\log Y) \otimes E, \nabla_{X/S}).$$

But it is an even stronger condition, permitting us to define a Gauss-Manin connection on these sheaves:

$$(3.4) \quad \mathrm{GM}(\nabla)^i : R^i f_*(\Omega^{\bullet}_{X/S}(\log Y) \otimes E, \nabla_{X/S})$$
$$\to \Omega^1_S(\log \Sigma) \otimes R^i f_*(\Omega^{\bullet}_{X/S}(\log Y) \otimes E, \nabla_{X/S}),$$

as the connecting morphism in relative cohomology of the exact sequence

$$(3.5) \quad 0 \to \Omega^1_S(\log \Sigma) \otimes \Omega^{\bullet-1}_{X/S}(\log Y) \otimes E \quad \to \quad \Omega^{\bullet}_X(\log Y)/ < f^* \Omega^2_S > \otimes E$$
$$\to \quad \Omega^{\bullet}_{X/S}(\log Y) \otimes E \to 0.$$

We want to prove formula (0.15)

$$\mathrm{N}w_n(\sum_i (-1)^i [R^i f_*(\Omega^{\bullet}_{X/S}(\log Y) \otimes E, \nabla_{X/S}), \mathrm{GM}^i(\nabla)])$$
$$= (-1)^d f_*(c_d(\Omega^1_{X/S}(\log Y/\Sigma), \mathrm{res}_{Z_i}) \cdot w_n((E, \nabla))) \in H^{2n}_{\mathrm{CS}}(S(\log(\Sigma)).$$

REDUCTION 3.2. *Assume $S = \mathrm{Spec}\,(F)$ for $F$ a function field.*

*Proof.* Equation 0.15 takes values in $H^{2p}_{\mathrm{CS}}(S(\log \Sigma))$. One applies Proposition 1.9. □

From now on, we will always assume that $S = \mathrm{Spec}\,(F)$ for $F$ a function field, even if strictly speaking it is a change of style, as $S$ is no longer a variety of finite type over $k$.

REDUCTION 3.3. *Assume now that $S = \mathrm{Spec}\,(F)$ is a field. Suppose further that $\dim X/S = d$ and that the Riemann-Roch theorem holds in dimensions $< d$. Then replace $E$ by $E(\sum m_i H_i)$ for $m_i \in \mathbb{Z}$, and $Y$ by $Y + \sum H_i$ where the $H_i$ are smooth divisors so that $Y + \sum H_i$ is a normal crossings divisor.*

*Proof.* By induction we may reduce to the case of a single very ample smooth divisor $H$. As in the proof of Proposition 1.4, we consider the tensor connection on $E(mH)$, still denoted by $\nabla$.

Consider first the case $m = 0$. The exact sequence of complexes

$$(3.6) \quad 0 \to \Omega^{\bullet}_{X/S}(\log Y) \otimes E \quad \to \quad \Omega^{\bullet}_{X/S}(\log(Y + H)) \otimes E$$
$$\to \quad \Omega^{\bullet}_{H/S}(\log(Y \cap H)) \otimes E[-1] \to 0$$

together with additivity of the Chern-Simons Newton class $\mathrm{N}w_*$ yields

$$(3.7) \quad \mathrm{N}w_n(\sum_i (-1)^i (R^i f_*(\Omega^{\bullet}_{X/S}(\log Y) \otimes E), \nabla_{X/S}), \mathrm{GM}^i(\nabla)))$$
$$= \mathrm{N}w_n(\sum_i (-1)^i (R^i f_*(\Omega^{\bullet}_{X/S}(\log(Y + H)) \otimes E), \nabla_{X/S}), \mathrm{GM}^i(\nabla)))$$
$$+ \mathrm{N}w_n(\sum_i (-1)^i (R^i f_*(\Omega^{\bullet}_{H/S}(\log(Y \cap H)) \otimes E), \nabla_{X/S}), \mathrm{GM}^i(\nabla))).$$

On the other hand, the image of the map

$$\Omega^1_{X/S}(\log Y)|_H \to \Omega^1_{X/S}(\log(Y+H))|_H$$

can be identified with $\Omega^1_{H/S}(\log(Y \cap H))$. Let $d = \dim X/S$, and let $i : H \to X$ be the inclusion. It follows from [27, Cor., p. 396], that in the group $\mathbb{H}^d(X, \underline{K}_{d,X,Y_\bullet})$,

$$(3.8)$$

$$c_d\Big(\Omega^1_{X/S}(\log(Y+H)), \mathrm{res}_Y\Big) = c_d\Big(\Omega^1_{X/S}(\log Y), \mathrm{res}_Y\Big)$$
$$+ i_* c_{d-1}\Big(\Omega^1_{H/S}(\log(Y \cap H)), \mathrm{res}_{H \cap Y}\Big).$$

Because the connection on $E$ is assumed regular along $H$, one has

$$(3.9) \quad c_d\Big(\Omega^1_{X/S}(\log(Y+H)), \mathrm{res}_{(Y+H)}\Big) \cdot \mathrm{N}w_n(E, \nabla)$$

$$= c_d\Big(\Omega^1_{X/S}(\log Y), \mathrm{res}_Y\Big) \cdot \mathrm{N}w_n(E, \nabla) \in H^{2(d+n)}_{\mathrm{CS}}(X(\log Y))$$

$$\subset H^{2(d+n)}_{\mathrm{CS}}(X(\log(Y+H))).$$

Also, with $g = f \circ i : H \to S$ and by our assumption that Riemann-Roch is true for fibre dimensions $< d$,

$$(3.10) \quad (-1)^{d-1} \mathrm{N}w_n\Big(\mathbb{R}g_*\Big(\Omega^*_{H/S}(\log(H \cap Y)) \otimes E\Big)\Big)$$

$$= g_*\Big(c_{d-1}\Big(\Omega^1_{H/S}(\log(Y \cap H)), \mathrm{res}_{H \cap Y}\Big) \cdot \mathrm{N}w_n((E, \nabla)|H)\Big)$$

$$= f_*\Big(i_* c_{d-1}\Big(\Omega^1_{H/S}(\log(Y \cap H)), \mathrm{res}_{H \cap D}\Big) \cdot \mathrm{N}w_n(E, \nabla)\Big).$$

Combining (3.7)–(3.10), we find that the Riemann-Roch formula (0.15) for $Y$ is equivalent to the Riemann-Roch formula (0.15) for $Y + H$.

We next consider the Riemann-Roch theorem for $(E(mH), \nabla)$ when $m \neq 0$. We show that formula (0.15) for $(E((m-1)H), \nabla)$ is equivalent to formula (0.15) for $(E(mH), \nabla)$. We claim first that

$$(3.11) \quad \mathrm{N}w_n(E, \nabla) = \mathrm{N}w_n(E(mH), \nabla)$$

$$\in H^{2n}_{\mathrm{CS}}(X(\log Y)) \subset H^{2n}_{\mathrm{CS}}(X(\log(Y+H))).$$

Indeed, it suffices to check this at the generic point. Since $\mathcal{O}_X(H)$ has a rational flat section, the bundles $E$ and $E(mH)$ are isomorphic (as bundles with connection) over $X - H$. It follows that the right-hand side of formula

(0.15) coincides for $E((m-1)H)$ and $E(mH)$. Next, there is an exact sequence of complexes
(3.12)

$$0 \rightarrow \Omega^*_{H/S}(\log(Y \cap H)) \otimes \mathcal{O}_H(H) \rightarrow \Omega^*_{X/S}(\log(Y+H))|_H \otimes \mathcal{O}_H(H)$$
$$\rightarrow \Omega^*_{H/S}(\log(Y \cap H)) \otimes \mathcal{O}_H(H)[-1] \rightarrow 0$$

which is compatible with the similar exact sequence where the forms relative to $S$ are replaced with absolute forms. It follows that

$$(3.13) \qquad \left[ \sum_i (-1)^i \{ R^i f_* \left( E(mH)/E((m-1)H) \otimes \Omega^*_{X/S}(\log(Y+H)), \nabla_{X/S}|_H \right), \mathrm{GM}^i(\nabla) \} \right] = 0 \in \mathcal{G}(S).$$

In particular the left-hand side of formula (0.15) coincides for $E((m-1)H)$ and $E(mH)$.                                  □

Recall that $E$ is locally free in formula (0.15).

REDUCTION 3.4. *Continue to assume* $S = \mathrm{Spec}\,(F)$ *for* $F$ *a field. Let* $\pi : X' \rightarrow X$ *be a birational morphism defined over* $F$, *such that* $X'$ *is smooth,* $Y' := \pi^{-1}(Y)$ *is a normal crossings divisor and* $\pi$ *is an isomorphism over* $U = X - Y$. *Then it suffices to prove* (0.15) *for* $(E', \nabla') = \pi^*(E, \nabla), f' = f \circ \pi$.

*Proof.* A well-known consequence of Deligne's mixed Hodge theory is that $R\pi_*\Omega^i_{X'/F}(\log Y') = \Omega^i_{X/F}(\log Y)$. In particular the projection formula applied to $\Omega^i_{X/F}(\log Y) \otimes \pi^*E$ enables one to identify

$$\mathbb{H}^r(X, \Omega^*_{X/S}(\log Y) \otimes E) \cong \mathbb{H}^r(X', \Omega^*_{X'/S}(\log Y') \otimes E')$$

as modules with connection on $S$, so the left-hand side of formula (0.15) for $(X, f, E, Y)$ and $(X', f', E', Y')$ is the same. A similar identification holds true for the right-hand side of Riemann-Roch. The pairing (2.10)

$$\Omega\underline{\underline{K}}_{n,X}(\log Y) \times \underline{\underline{K}}_{d,X,Y_\bullet} \rightarrow \Omega\underline{\underline{K}}_{d+n,X}$$

maps via $R\pi_*$ to the pairing

$$R\pi_*\Omega\underline{\underline{K}}_{n,X'}(\log Y') \times R\pi_*\underline{\underline{K}}_{l,X',Y'_\bullet} \rightarrow R\pi_*\Omega\underline{\underline{K}}_{d+n,X'},$$

and

$$(3.14) \qquad \Gamma(X, \underline{\underline{H}}^{d+n}(\Omega\underline{\underline{K}}_{d+n,X})) \quad \subset \quad \Gamma(X', \underline{\underline{H}}^{d+n}(\Omega\underline{\underline{K}}_{d+n,X'}))$$
$$\subset \quad \Gamma(U, \underline{\underline{H}}^{d+n}(\Omega\underline{\underline{K}}_{d+n,X}))$$

by Proposition 1.9. On the other hand, the class $\mathrm{N}w_n(E, \nabla)$ comes from $\mathbb{H}^n(X, \Omega\underline{\underline{K}}_{n,X}(\log Y))$; thus the class $\mathrm{N}w_n(E', \nabla').c_d(\Omega^1_{X'/F}(\log Y'))$ comes from $\mathbb{H}^{n+d}(X', \Omega\underline{\underline{K}}_{n+d,X'})$. Denoting by $j : U \rightarrow X$ and $j' : U \rightarrow X'$ the open embeddings, one has an exact triangle

$$\Omega\underline{\underline{K}}_{n+d,X} \rightarrow R\pi_*\Omega\underline{\underline{K}}_{n+d,X'} \rightarrow R\pi_*j'_*\underline{\underline{K}}_{n+d}/j_*\underline{\underline{K}}_{n+d}$$

and since $n + d > d$, one has $\mathbb{H}^{d+n}(V, R\pi_* j'_* \underline{\underline{K}}_{n+d} / j_* \underline{\underline{K}}_{n+d}) = 0$ for any Zariski open set $V \subset X$. This shows that

$$(3.15) \qquad \mathrm{Im}\ \mathbb{H}^{n+d}(X, \Omega\underline{\underline{K}}_{n+d,X}) \quad = \quad \mathrm{Im}\ \mathbb{H}^{n+d}(X', \Omega\underline{\underline{K}}_{n+d,X'})$$
$$\mathrm{in}\ \mathbb{H}^{d+n}(U, \Omega\underline{\underline{K}}_{d+n,X}).$$

Since on $U$, one trivially has

$$c_d(\Omega^1_{X'/F}(\log Y'), \mathrm{res}) = c_d(\Omega^1_{X/F}(\log Y), \mathrm{res}),$$

where the maps res on $X$ (resp. $X'$) take into account all the components of $Y$ (resp. $Y'$), one concludes that

$$(3.16) \quad \mathrm{N}w_n(E, \nabla).c_d(\Omega^1_{X/F}(\log Y)) \quad = \quad \mathrm{N}w_n(E', \nabla').c_d(\Omega^1_{X'/F}(\log Y'))$$
$$\in\ \Gamma(X, \underline{\underline{H}}^{d+n}(\Omega\underline{\underline{K}}_{d+n,X})).$$

This shows that the right-hand side of formula $(0.15)$ is the same for $E$ and $E'$. $\qquad \square$

REDUCTION 3.5. *With notation as above, it suffices to prove formula* $(0.15)$ *in the case* $S = \mathrm{Spec}\,(F)$ *for* $F$ *a field, and* $\dim X/S = 1$.

*Proof.* We have already reduced to the case $S = \mathrm{Spec}\,(F)$. Assume $d = \dim X/S > 1$. Using Reduction 3.4, we may blow up the base of a Lefschetz pencil, and assume we have a factorization

$$(3.17) \qquad\qquad\qquad X \xrightarrow{g} \mathbb{P}^1_S \xrightarrow{h} S$$

with $f = h \circ g$. By Lemma 2.3, the sheaf $\Omega^1_{X/\mathbb{P}^1_S}(\log Y)$ is locally free away from $Y_{i_1} \cap \ldots \cap Y_{i_d}$, of codimension $d$, finite over $S$, and away from the singularities of the morphisms $g_i : Y_i \to \mathbb{P}^1_S$. By Reduction 3.4, we may blow up the intersections $Y_{i_1} \cap \ldots \cap Y_{i_d}$ and replace Y by its total transform. We further blow up the singularities of the bad fibers of $g_i$ in $X$, so that the total inverse image of $Y_i$ becomes a normal crossings divisor. Again by Reduction 3.4, we may replace $X$ by this blowup and $Y_i$ by its total transform. Thus we may assume that

$$\{g : X \to \mathbb{P}^1_S, Y, \Sigma\}$$

is a relative normal crossings divisor. We write $Y = g^{-1}(\Sigma)_{\mathrm{red}} + Z$. We next have to show that the curvature condition $(3.1)$ is fulfilled for the morphism $g$. Let $S^2$ be the 2-nd symmetric tensor power of $\Omega^1_{\mathbb{P}^1_S}(\log \Sigma)$. In order to simplify the notation, we set $\Omega^i_{\mathbb{P}^1_S}(\log \Sigma) = \Omega^i_b$, $\Omega^i_X(\log Y) = \Omega^i_s$, $\Omega^i_{X/\mathbb{P}^1_S}(\log Y) = \Omega^i_r$ in the following commutative diagram:

(3.18)

$$
\begin{array}{ccc}
(\Omega_b^2/\Omega_F^2) \otimes \Omega_r^{*-2} \otimes E & = & (\Omega_b^2/\Omega_F^2) \otimes \Omega_r^{*-2} \otimes E \\
\downarrow & & \downarrow \\
0 \to \left((\Omega_b^1 \otimes \Omega_s^{*-1}/<\Omega_b^2>)/<\Omega_F^2, S^2>\right) \otimes E & \to (\Omega_s^*/<\Omega_F^2, \Omega_b^3>) \otimes E & \to \Omega_r^* \otimes E \to 0 \\
\downarrow & & \downarrow \quad\quad \downarrow \\
0 \to \Omega_b^1 \otimes \Omega_r^{*-1} \otimes E & \to \Omega_s^*/<\Omega_b^2> \otimes E & \to \Omega_r^* \otimes E \to 0 \,.
\end{array}
$$

This shows that the composite of the connecting morphisms

$$
(3.19) \qquad
\begin{aligned}
R^i f_*(\Omega_r^* \otimes E) & \to R^i f_*(\Omega_b^1 \otimes \Omega_r^{*-1} \otimes E) \\
& \to R^i f_*((\Omega_b^2/\Omega_F^2) \otimes \Omega_r^{*-2} \otimes E)
\end{aligned}
$$

is vanishing.

By induction on $d$ we may assume the Riemann-Roch Theorem 0.1 with values in

$$
H_{\mathrm{CS}}^{2n}(\mathbb{P}_S^1(\log \Sigma))
$$

holds for $(E, \nabla)$ and the morphism $g$.

LEMMA 3.6. *Let $\mathcal{F}$ be a coherent sheaf on $\mathbb{P}_F^1$ and let*

$$
\nabla : \mathcal{F} \to \mathcal{F} \otimes \Omega_{\mathbb{P}_F^1}^1(\log \Sigma)
$$

*be a connection, and assume $\nabla^2(\mathcal{F}) \subset \mathcal{F} \otimes h^*\Omega_F^2$. Then the torsion subsheaf $\mathcal{F}_{\mathrm{tors}} \subset \mathcal{F}$ is stable under $\nabla$, and*

$$
(3.20) \qquad \mathbb{R}^0 h_*(\mathcal{F}_{\mathrm{tors}} \otimes \Omega_{\mathbb{P}_F^1/F}^*(\log \Sigma)) - \mathbb{R}^1 h_*(\mathcal{F}_{\mathrm{tors}} \otimes \Omega_{\mathbb{P}_F^1/F}^*(\log \Sigma)) = 0
$$

*in the Grothendieck group of finite dimensional $F$-vector spaces with connection.*

*Proof.* Note first that, as is well-known, the support of $\mathcal{F}_{\mathrm{tors}}$ is contained in $\Sigma$. Indeed, if $t$ is a local parameter at a point not in $\Sigma$ and $t^n \mathcal{F}_{\mathrm{tors}} = (0)$ for some $n > 0$, we have for $s$ a torsion section,

$$
0 = \nabla(t^n s) = n t^{n-1} s \otimes dt + t^n \nabla(s).
$$

Multiplying through by $t$, we see that $t^{n+1}\nabla(s) = 0$, so that $\nabla(s)$ is torsion, and $t^n \nabla(s) = 0$. It follows that $t^{n-1}\mathcal{F}_{\mathrm{tors}} = (0)$.

Now suppose $t$ is a local parameter at a point of $\Sigma$. Replacing $dt$ by $dt/t$ in the above equation, we see that $\nabla(\mathcal{F}_{\mathrm{tors}}) \subset \mathcal{F}_{\mathrm{tors}}$, and $\nabla$ stabilizes the filtration

$$
N^i \mathcal{F}_{\mathrm{tors}} = \{\varphi \in \mathcal{F}_{\mathrm{tors}}, t^i \varphi = 0\}.
$$

One is thus reduced to showing (3.20) in the case when $\Sigma$ is a single closed point and $\mathcal{F}_{\mathrm{tors}}$ is an $F(\Sigma)$-vector space. We will use the residue along $\Sigma$ to split the Gauss-Manin on $\mathcal{F}_{\mathrm{tors}}$. Write $M := \mathcal{F}_{\mathrm{tors}}$ and $L := F(\Sigma)$. We have

$$M \otimes \Omega^1_{\mathbb{P}^1_F}(\log \Sigma) \cong (M \otimes_L \Omega^1_L) \oplus M,$$

where projection onto the second factor on the right corresponds to taking the residue at $\Sigma$, and the splitting depends on the choice of $t$. The absolute connection is given by a pair $(A, B)$ with $A : M \to M \otimes \Omega^1_L$, and $B : M \to M$. (The curvature condition means $(B \otimes 1)A = AB$.) To calculate the Gauss-Manin connection, note $M \otimes \Omega^1_{\mathbb{P}^1_F/F}(\log \Sigma) \cong M$. The exact sequence of absolute to relative differentials, coupled to $M$, yields a diagram (with $\sigma$ being the evident splitting)

$$
\begin{array}{ccc}
M & =\!=\!= & M \\
\downarrow{\scriptstyle A \oplus B} & & \downarrow{\scriptstyle B} \\
\end{array}
$$

$$
\begin{array}{ccccc}
M \otimes_F \Omega^1_F & \longrightarrow & (M \otimes \Omega^1_L) \oplus M & \xrightarrow{\ \sigma\ } & M \\
\downarrow{\scriptstyle B} & & \downarrow{\scriptstyle (-B \otimes 1, A)} & & \\
M \otimes \Omega^1_L & \xrightarrow{\ \cong\ } & M \otimes \Omega^1_L. & & \\
\end{array}
$$

Viewing this as an exact sequence of complexes and taking boundaries, we find a representative for $[\mathbb{R}h_*(M \otimes \Omega^*_{\mathbb{P}^1_F/F}(\log(\Sigma))]$ of the form

$$
\begin{array}{ccc}
M & \xrightarrow{\ A\ } & M \otimes_F \Omega^1_F \\
\downarrow{\scriptstyle B} & & \downarrow{\scriptstyle B \otimes 1} \\
M & \xrightarrow{\ A\ } & M \otimes_F \Omega^1_F. \\
\end{array}
$$

Since the top and bottom rows are the same connection, the total class in the Grothendieck group is zero.                                                      $\square$

Define

$$\mathcal{F}^{\mathrm{ev}} := \left( \oplus_{i \geq 0} \mathbb{R}^{2i} g_*(E \otimes \Omega^*_{X/\mathbb{P}^1_S}(\log Y)) \right.$$

$$\left. \oplus_{i \geq 0} \mathbb{R}^{2i+1} g_*(\mathcal{O}_X^{\oplus \mathrm{rk}(E)} \otimes \Omega^*_{X/\mathbb{P}^1_S}(\log Y)) \right) \Big/ (\mathrm{torsion}),$$

$$\mathcal{F}^{\mathrm{odd}} := \left( \oplus_{i \geq 0} \mathbb{R}^{2i+1} g_*(E \otimes \Omega^*_{X/\mathbb{P}^1_S}(\log Y)) \right.$$

$$\left. \oplus_{i \geq 0} \mathbb{R}^{2i} g_*(\mathcal{O}_X^{\oplus \mathrm{rk}(E)} \otimes \Omega^*_{X/\mathbb{P}^1_S}(\log Y)) \right) \Big/ (\mathrm{torsion}).$$

Applying the Hirzebruch Riemann-Roch theorem at the generic point of $\mathbb{P}^1_S$ (compare (0.3)) implies $\mathrm{rk}(\mathcal{F}^{\mathrm{ev}}) = \mathrm{rk}(\mathcal{F}^{\mathrm{odd}})$. By the lemma, we have, with $\tilde{E} = E - \mathrm{rk}(E)\mathcal{O}_{\mathbb{P}^1_S}$,

$$\left[ \mathbb{R}g_*(\Omega^*_{X/\mathbb{P}^1_S}(\log Y) \otimes \tilde{E}), \mathrm{GM}(\nabla) \right] = [\mathcal{F}^{\mathrm{ev}}, \nabla^{\mathrm{ev}}] - [\mathcal{F}^{\mathrm{odd}}, \nabla^{\mathrm{odd}}].$$

The sheaves $\mathcal{F}^{\mathrm{ev}}, \mathcal{F}^{\mathrm{odd}}$ are locally free and their connections satisfy the basic curvature form condition by (3.19). Hence we have

(3.21)

$$\mathrm{Nw}_n\left( \mathbb{R}f_*(\Omega^*_{X/S}(\log Y) \otimes \tilde{E}), \mathrm{GM}(\nabla) \right)$$

$$= \mathrm{Nw}_n\left( \mathbb{R}h_*(\Omega^*_{\mathbb{P}^1_S/S}(\log \Sigma) \otimes \mathcal{F}^{\mathrm{ev}}) \right) - \mathrm{Nw}_n\left( \mathbb{R}h_*(\Omega^*_{\mathbb{P}^1_S/S}(\log \Sigma) \otimes \mathcal{F}^{\mathrm{odd}}) \right)$$

$$= h_*\left( (\mathrm{Nw}_p(\mathcal{F}^{\mathrm{ev}}, \nabla^{\mathrm{ev}}) - \mathrm{Nw}_p(\mathcal{F}^{\mathrm{odd}}, \nabla^{\mathrm{odd}})) \cdot c_1(\Omega^1_{\mathbb{P}^1_S/S}(\log \Sigma), \mathrm{res}_\Sigma) \right)$$

$$= h_*\left( g_*\left( \mathrm{Nw}_p(E, \nabla) \cdot c_{d-1}(\Omega^1_{X/\mathbb{P}^1_S}(\log Y), \mathrm{res}_Z) \right) \cdot c_1(\Omega^1_{\mathbb{P}^1_S/S}(\log \Sigma), \mathrm{res}_\Sigma) \right)$$

$$= h_*g_*\left( \mathrm{Nw}_p(E, \nabla) \cdot c_{d-1}(\Omega^1_{X/\mathbb{P}^1_S}(\log Y), \mathrm{res}_Z) \cdot g^*c_1(\Omega^1_{\mathbb{P}^1_S/S}(\log \Sigma), \mathrm{res}_\Sigma) \right)$$

$$\overset{2.16, 2.8}{=} f_*\left( \mathrm{Nw}_p(E, \nabla) \cdot c_d(\Omega^1_{X/S}(\log Y), \mathrm{res}_Y) \right). \qquad \square$$

REDUCTION 3.7. *It suffices to prove formula* (0.15) *in the case*

$$S = \mathrm{Spec}\,(F), \ X = \mathbb{P}^1_F.$$

*Proof.* By Reduction 3.5, it suffices to prove the formula for $f : X \to \mathrm{Spec}\,(F)$ a complete smooth curve. We factor $f$ as follows

$$X \overset{g}{\to} \mathbb{P}^1_F \overset{h}{\to} \mathrm{Spec}\,(F)$$

where $g$ is finite. Enlarging $Y$ if necessary, we get for a suitable subscheme $\Sigma \subset \mathbb{P}^1_F$, finite over $F$,

$$g^*(\Omega^1_{\mathbb{P}^1_F/F}(\log \Sigma)) \cong \Omega^1_{X/F}(\log Y).$$

Projection formulae give

(3.22) $$\mathbb{R}f_*(\Omega^*_{X/F}(\log Y) \otimes \tilde{E}) = \mathbb{R}h_*(\Omega^*_{\mathbb{P}^1_F/F}(\log \Sigma) \otimes g_*\tilde{E})$$

(3.23) $$f_*\left( \mathrm{Nw}_n(E, \nabla) \cdot c_1(\Omega^1_{X/F}(\log Y), \mathrm{res}_Y) \right)$$

$$= h_*\left( g_*(\mathrm{Nw}_n(E, \nabla)) \cdot c_1(\Omega^1_{\mathbb{P}^1_F/F}(\log \Sigma), \mathrm{res}_\Sigma) \right).$$

We have reduced the problem to showing

$$(3.24) \qquad g_*(\mathrm{N}w_n(E, \nabla)) = \mathrm{N}w_n(g_*\tilde{E}, g_*\nabla).$$

Since our classes are recognized at the generic point of the variety, we reduce to the case where

$$g : \mathrm{Spec}\,(M) \to \mathrm{Spec}\,(L)$$

is a finite map, where $L$ is a function field over $F$, $M = L[t]/\langle\varphi(t)\rangle$ is a commutative, semi-simple $L$-algebra, $\varphi(t)$ is a polynomial of degree $r \geq 1$, $E = \oplus_1^N M$, with basis $e_i$, and $\nabla$ is given by an $N \times N$ matrix $A(t) = \sum_{i=0}^{i=r-1} t^i A_i$. Let $L' = L(a_1, \dots, a_r)$ be the Galois hull of $M$, with $\varphi(a_i) = 0$, $M' = L' \otimes_L M = L'[t]/\langle\varphi(t)\rangle = \Pi_{j=1}^{j=r} L'_j$ where the projection on the $j^{\mathrm{th}}$ factor $L'_j \cong L'$ is induced by $t \mapsto a_j$. As the receiving group for $\mathrm{N}w_n, n > 1$, is torsion-free, and both terms of formula (3.24) are compatible with base-change, we are reduced to showing the formula for $g' : \mathrm{Spec}\,(M') \to \mathrm{Spec}\,(L')$, $(E', \nabla')$, with $(E', \nabla')|\mathrm{Spec}\,(L'_j) = (\oplus_1^N L', A(a_j))$. Then the left-hand side of formula (3.24) becomes

$$\sum_{j=1}^{j=r} w_n((E', \nabla')|\mathrm{Spec}\,(L'_j)) = \sum_{j=1}^{j=r} w_n(A(a_j)),$$

whereas

$$(g'_* E', g_*\nabla) = (\sum_1^{rN} L', \mathrm{diag}(A(a_1), \dots, A(a_r))),$$

and thus the right-hand side is $\sum_{j=1}^{j=r} w_n(A(a_j))$. This concludes the proof for $n > 1$, and for $n = 1$ as well, but modulo torsion.

In order to understand the torsion-factor, one does the following direct calculation. Let $\alpha_i \in L$ be the trace of the matrix $A_i$.

$$\mathrm{Tr}_{M/L}(\sum_{i=0}^{i=r-1} t^i \alpha_i).$$

On the other hand, consider the $L$-basis

$$e_1, \dots, e_N, te_1, \dots, te_N, \dots, t^{r-1}e_1, \dots, t^{r-1}e_N$$

of $g_*E$. In this basis, $g_*\nabla$ is an $r \times r$ block matrix, each block being of size $N \times N$. The Leibniz formula applied to $t^i e_j$ implies then that

$$w_1(g_*E, g_*\nabla) = \mathrm{Tr}_{M/L}(\sum_{i=0}^{i=r-1} t^i \alpha_i) + N \sum_{j=0}^{j=r-1} \beta_{jj},$$

where $d(t^j) = \sum_{i=0}^{i=r-1} t^i \beta_{ji}$. In other words, one obtains the formula

$$(3.25) \qquad g_*(w_1(E, \nabla)) + \mathrm{N}w_1(g_*M, g_*d) = w_1(g_*E, g_*\nabla).$$

It remains to show that $w_1(g_*M, g_*d)$ is 2-torsion. Since $d_L \circ \text{Tr}_{M/L} = \text{Tr}_{M/L} \circ d_M$, the pairing $M \otimes_L M \to L, (x, y) \mapsto \text{Tr}_{M/L}(xy)$ induces an isomorphism between $(g_*M, g_*d)$ and its dual. With $B$ for the connection matrix of $g_*M$ in some basis, the connection matrix of the dual in the dual basis is $-B^t$. Thus for some invertible matrix $\phi$ with coefficients in $L$ we get

$$-B^t = \phi B \phi^{-1} + d\phi \cdot \phi^{-1}.$$

Taking traces,

$$2w_1(g_*M, g_*\nabla) = -\text{Tr}(d\phi \cdot \phi^{-1}) = -d\log(\det(\phi)) \mapsto 0 \in H^2_{\text{CS}}(\text{Spec}(L)). \quad \square$$

## 4. The Riemann-Roch theorem for $\mathbb{P}^1$

In this section $F \supset k$ is a field, $P := \mathbb{P}^1_F \xrightarrow{f} \text{Spec}(F)$, and $D \subset P$ is a reduced, effective divisor. We are given $E$ a vector bundle of rank $N$ on $P$ with a $k$-connection

$$\nabla : E \to E \otimes \Omega^1_{P/k}(\log D),$$

with curvature

$$(4.1) \qquad\qquad \nabla^2 \in f^*\Omega^2_F \otimes \text{End}(E).$$

We wish to prove the Riemann-Roch Theorem 0.1 for $(E, \nabla)$. By Reduction 3.3 we may tensor by a multiple of $\mathcal{O}_P(1 \cdot \infty)$ and assume

$$(4.2) \qquad E \cong \mathcal{O}(m_1) \oplus \ldots \oplus \mathcal{O}(m_r); \ 0 = m_1 \leq m_2 \leq \ldots \leq m_r.$$

REDUCTION 4.1. *It suffices to prove the Riemann-Roch theorem for* $E = \oplus_r \mathcal{O}_P$.

*Proof.* Changing coordinates if necessary, we may assume $\infty \notin D$. By Reduction 3.3, we can replace $D$ by $D + \infty$, so that $\infty \in D$ and $\nabla$ has trivial residue at $\infty$. We think of $E$ as a direct sum of line bundles associated to divisors supported at infinity

$$E \cong \oplus \mathcal{O}(m_i \cdot \infty).$$

Define

$$E' = \bigoplus_{m_i = 0} \mathcal{O}_P \oplus \bigoplus_{m_i > 0} \mathcal{O}((m_i - 1) \cdot \infty) \subset E.$$

We shall show that triviality of the residue at $\infty$, together with the basic curvature form hypothesis Definition 3.1, implies that $E'$ is stable under the connection. Since we have normalized so that $m_1 = 0$ and all $m_i \geq 0$, the number $m_r$ is well-defined. By induction on $m_r$, Riemann-Roch holds for $E'$. But Lemma 3.6 implies that the left-hand side of the Riemann-Roch formula (0.15) coincides for $E$ and $E'$. The same holds for the right-hand side because $E \cong E'$ away from $\infty$, so that Riemann-Roch holds for $E$.

LEMMA 4.2.   *With notation as above, $E' \subset E$ is stable under $\nabla$.*

*Proof.* The assertion is invariant under an extension of $F$, so we may assume $D = \{a_1, \ldots, a_\delta, \infty\}$ with all $a_\nu \in P(F)$. Let $1_j$, $1 \leq j \leq r$, be the evident basis of $E$ on $\mathbb{A}^1 = P - \{\infty\}$, and let $z$ be the standard parameter on $P$. An element $\gamma \in \Gamma(P, \mathcal{O}(n) \otimes \Omega^1_{\mathbb{P}^1_F}(\log D))$ for $n \geq 0$ can be uniquely written in the form

$$\sum_{\nu=1}^{\delta} A_\nu d\log(z - a_\nu) + \sum_{i=0}^{n} z^i \eta_i + \sum_{j=1}^{n} C_j (z - a_1)^j d\log(z - a_1),$$

with $A_\nu, C_j \in F$ and $\eta_i \in \Omega^1_F$. Since $\nabla(z^{m_j} 1_j) \in \Gamma(E \otimes \Omega^1_{\mathbb{P}^1_F}(\log D))$, we may write

$$
\begin{aligned}
(4.3) \quad \nabla(1_j) \;\; = \;\; & \sum 1_k \otimes \Big[ \sum_{i=0}^{m_k - m_j} \eta_j^{ki}(z - a_1)^i + \sum A_j^{k\nu} d\log(z - a_\nu) \\
& \qquad + \sum_{\ell=1}^{m_k - m_j} C_j^{k\ell}(z - a_1)^\ell d\log(z - a_1) \Big].
\end{aligned}
$$

If $m_j > m_k$, the sums over $i$ and $\ell$ on the right are not there. If $m_j = m_k$, the sum over $\ell$ is absent. With respect to (4.3) we have the following facts:

$$(4.4) \qquad\qquad\qquad C_j^{k, m_k - m_j} = 0.$$

For $m_j \geq m_k$,

$$(4.5) \qquad \sum_\nu A_j^{k\nu} d\log(z - a_\nu) \in \Gamma(\mathbb{P}^1, \Omega^1_{\mathbb{P}^1_F}(\log D)((m_k - m_j) \cdot \infty)).$$

$$(4.6) \qquad \sum A_j^{k\nu} d\log(z - a_\nu) \in \Gamma(\mathbb{P}^1, \Omega^1_{\mathbb{P}^1_F/F}(\log D)((m_k - m_j - 1) \cdot \infty)).$$

To check (4.4) we may suppose $m_k > m_j$. The composition

$$
\begin{aligned}
(4.7) \quad \mathcal{O}(m_j) \quad &\hookrightarrow \quad E \xrightarrow{\nabla} E \otimes \Omega^1_{\mathbb{P}^1_F}(\log D) \to E \otimes \Omega^1_{\mathbb{P}^1_F/F}(\log D) \\
&\xrightarrow{\mathrm{res}_\infty} \quad E|_\infty \twoheadrightarrow \mathcal{O}(m_k)|_\infty = \mathcal{O}(m_k \cdot \infty)/\mathcal{O}((m_k - 1) \cdot \infty)
\end{aligned}
$$

maps

$$(z - a_1)^{m_j} 1_j \mapsto C_j^{k, m_k - m_j}(z - a_1)^{m_k} \pmod{\mathcal{O}((m_k - 1) \cdot \infty)}.$$

By assumption, the connection has zero residue at infinity, so this is zero. The inclusion (4.5) follows because

$$
\begin{aligned}
\nabla((z - a_1)^{m_j} 1_j) \;\; = \;\; & \sum (z - a_1)^{m_k} 1_k \\
& \otimes \Big[ \delta_{m_j, m_k} \cdot \eta_j^{k0} + (z - a_1)^{m_j - m_k} \sum A_j^{k\nu} d\log(z - a_\nu) \Big] \\
& + m_j(z - a_1)^{m_j} 1_j \otimes d\log(z - a_1)
\end{aligned}
$$

is assumed to extend across infinity. Finally, (4.6) holds because of the vanishing of the residue (4.7). In the case $m_j \geq m_k$ the residue map is

$$1_j \mapsto \left( (z - a_1)^{m_j - m_k} \sum A_j^{k\nu} d\log(z - a_\nu) \right)\Big|_\infty .$$

Now view (4.3) as defining the connection matrix $B = (b_j^k)$ for $\nabla$ on $P - \{\infty\}$. The above assertions can be summarized as follows. For $m_k > m_j$,

$$b_j^k \in \Gamma\left( \mathbb{P}^1, f^*\Omega_F^1((m_k - m_j)\cdot\infty) + \Omega_{\mathbb{P}_F^1}^1(\log D)((m_k - m_j - 1)\cdot\infty) \right)$$

and for $\eta \in \Omega_F^1$ and $m_k \leq m_j$,

$$b_j^k \in \Gamma\left( \mathbb{P}^1, \Omega_{\mathbb{P}_F^1}^1(\log D)((m_k - m_j)\cdot\infty) \right),$$

$$b_j^k \in \Gamma\left( \mathbb{P}^1, \Omega_{\mathbb{P}_F^1/F}^1(\log D)((m_k - m_j - 1)\cdot\infty) \right),$$

$$\eta \wedge b_j^k \bmod f^*\Omega_F^2 \in \Gamma\left( \mathbb{P}^1, f^*\Omega_F^1 \otimes \Omega_{\mathbb{P}_F^1/F}^1(\log D)((m_k - m_j - 1)\cdot\infty) \right).$$

It follows that whenever $m_i \leq m_k$,

$$b_i^j b_j^k \bmod f^*\Omega_F^2 \in \Gamma(\mathbb{P}^1, f^*\Omega_F^1 \otimes \Omega_{\mathbb{P}_F^1/F}^1(\log D)(m_k - m_i - 1)\cdot\infty).$$

Vanishing of $C_i^{k,m_k - m_i}$ and the basic curvature form condition

$$db_i^k \equiv \sum_j b_i^j b_j^k \bmod (f^*\Omega_F^2)$$

imply

$$db_i^k = (m_k - m_i)\eta_i^{k,m_k - m_i}(z - a_1)^{m_k - m_i} d\log(z - a_1) + \varepsilon$$

for $\varepsilon \in \Gamma\left( \mathbb{P}^1, (\Omega_{\mathbb{P}_F^1}^2(\log(D))/f^*\Omega_F^2)((m_k - m_i - 1)\infty) \right)$. It follows for $m_k > m_i$ that $\eta_i^{k,m_k - m_i} = 0$ as claimed.

It follows now from (4.3) that $\nabla$ stabilizes $E' \subset E$, proving the lemma. $\quad\square$

As mentioned, the lemma implies Reduction 4.1 by induction on $m_r$. $\quad\square$

We assume now $E = \mathcal{O}_P^{\oplus N}$ with connection given by

$$(4.8) \qquad \nabla(1_j) = \sum_{k\nu} A_j^{k\nu} 1_k \otimes d\log(z - a_\nu) + \sum_k 1_k \otimes \eta_j^k,$$

where $a_i \in \mathbb{P}^1(S)$. Note that this assumption of $S$-rationality of the $a_i$ implies that the correcting term $\operatorname{rank}(E) \cdot \mathrm{Nw}_p(H_{\mathrm{DR}}^*(\mathbb{P}_S^1/S(\log D)), \nabla_{\mathrm{GM}})$ of formula (0.15) vanishes, as it trivially vanishes for $D = \emptyset$, and $f|Y*(\mathcal{O}_D, d_D)$ is the trivial connection of rank $= \deg D$. We will remove this hypothesis later (cf. Lemma 4.8).

Our objective is to compute the Gauss-Manin connection. First, we define an $F$-linear splitting $\sigma$ (already used above) of the natural reduction from absolute to relative $E$-valued 1-forms

$$
\begin{aligned}
(4.9) \qquad \sigma : \Gamma(P, \Omega^1_{P/F}(\log D)) \quad &\to \quad \Gamma(P, \Omega^1_{P/k}(\log D)); \\
\sigma(1_k \otimes (z - a_1)^\nu dz) \quad &= \quad 1_k \otimes (z - a_1)^\nu d(z - a_1), \\
\sigma(1_k \otimes d\log(z - a_\nu)) \quad &= \quad 1_k \otimes d\log(z - a_\nu),
\end{aligned}
$$

Now consider the diagram
(4.10)

$$
\begin{array}{ccc}
\Gamma(E) & = & \Gamma(E) \\
\downarrow \nabla_1 & & \downarrow \nabla_{P/F} \\
\Gamma(E) \otimes \Omega^1_F \quad \to \quad \Gamma(E \otimes \Omega^1_P(\log D)) & \xrightarrow{\sigma} & \Gamma(E \otimes \Omega^1_{P/F}(\log D)) \\
\downarrow \nabla_{P/F} \otimes 1 \qquad\qquad \downarrow \nabla_2 & & \\
\Gamma(E \otimes \Omega^1_{P/F}(\log D)) \otimes \Omega^1_F & = & \Gamma(E \otimes \Omega^2_P(\log D)/\Omega^2_F).
\end{array}
$$

Here $\nabla_1$ and $\nabla_2$ are the absolute connection maps. Define

$$
(4.11) \qquad \Phi := \nabla_1 - \sigma \nabla_{P/F}; \quad \Psi = -\nabla_2 \sigma.
$$

The diagram

$$
(4.12) \qquad
\begin{array}{ccc}
\Gamma(E) & \xrightarrow{\Phi} & \Gamma(E) \otimes \Omega^1_F \\
\nabla_{P/F} \downarrow & & \downarrow \nabla_{P/F} \otimes 1 \\
\Gamma(E \otimes \Omega^1_{P/F}(\log D)) & \xrightarrow{\Psi} & \Gamma(E \otimes \Omega^1_{P/F}(\log D)) \otimes \Omega^1_F
\end{array}
$$

represents $\sum_{i=0}^{i=2} (-1)^i (R^i f_*(E \otimes \Omega^*_{P/F}, \nabla_{P/F}), \mathrm{GM}^i(\nabla))$ in the Grothendieck group $\mathcal{K}(F)$ of $F$-vector spaces with connection. We see from (4.8) that

$$
(4.13) \qquad \Phi(1_j) = (\nabla - \sigma \nabla_{P/F})(1_j) = \sum_k 1_k \otimes \eta^k_j.
$$

Also

$$
\begin{aligned}
(4.14) \quad \Psi(1_j \otimes d\log(z - a_\nu)) &= -\nabla(1_j) \wedge d\log(z - a_\nu) \\
&= \left( -\sum_{k\tau} A^{k\tau}_j 1_k \otimes d\log(z - a_\tau) - 1_k \otimes \eta^k_j \right) \wedge d\log(z - a_\nu) \\
&= -\sum_{k\tau} A^{k\tau}_j 1_k \otimes \left( d\log(z - a_\tau) - d\log(z - a_\nu) \right) \otimes d\log(a_\nu - a_\tau) \\
&\quad + \sum_k 1_k \otimes d\log(z - a_\nu) \otimes \eta^k_j.
\end{aligned}
$$

Define

$$(4.15) \qquad B_{j\nu}^{k\tau} = \begin{cases} -A_j^{k\tau} d\log(a_\nu - a_\tau) & \tau \neq \nu \\ \eta_j^k + \sum_{\theta \neq \nu} A_j^{k\theta} d\log(a_\nu - a_\theta) & \tau = \nu \end{cases}.$$

Then

$$(4.16) \qquad \Psi(1_j \otimes d\log(z - a_\nu)) = \sum_{k\tau} 1_k \otimes d\log(z - a_\tau) \wedge B_{j\nu}^{k\tau}.$$

The left-hand side of formula (0.15) for $E$ (in this case there is no need to subtract off $\mathcal{O}^{\oplus \mathrm{rk}(E)}$) is given by

$$(4.17)$$
$$\mathrm{Nw}_n((\sum_{i=0}^{i=2}(-1)^i(R^i f_*(E \otimes \Omega_{P/F}^*, \nabla_{P/F}), \mathrm{GM}^i(\nabla))) = \mathrm{Nw}_n(\Phi) - \mathrm{Nw}_n(\Psi).$$

We next make some observations about $\mathrm{Nw}_n(\Psi)$. Define $B_\nu^\tau$ (resp. $B$) to be the $N \times N$ matrix (resp. $\delta \times \delta$ block matrix with blocks of size $N \times N$)

$$(4.18) \qquad B_\nu^\tau := (B_{j\nu}^{k\tau})_{1 \leq j,k \leq N} \quad (\mathrm{resp.} B = (B_\nu^\tau)_{1 \leq \nu,\tau \leq n}).$$

LEMMA 4.3. *Let* $M(B) = B^{r_1}(dB)^{r_2} \cdots B^{r_{2s-1}}(dB)^{r_{2s}}$ *be some* (*noncommuting*) *monomial in* $B$ *and* $dB$. *Then*

$$\mathrm{Tr} \ (M(B)) = \sum_{\tau=1}^n \mathrm{Tr} \ (M(B_\tau^\tau)).$$

*Proof.* Write as above

$$M(B)_\nu^\tau := (M(B)_{j\nu}^{k\tau})_{1 \leq j,k \leq N} \quad (\mathrm{resp.} M(B) = (M(B)_\nu^\tau)_{1 \leq \nu,\tau \leq n}).$$

Then $\mathrm{Trace}(M(B)) = \sum_\tau \mathrm{Trace}(M(B)_\tau^\tau)$. Now

$$M(B)_\tau^\tau = \sum_{\tau_1, \ldots, \tau_{r_{2s}-1}} B_\tau^{\tau_1} B_{\tau_1}^{\tau_2} \cdots B_{\tau_{r_1-1}}^{\tau_{r_1}} dB_{\tau_{r_1}}^{\tau_{r_1+1}} \cdots dB_{\tau_{r_{2s}-1}}^\tau.$$

For $\nu \neq \tau$ we can write $B_\nu^\tau = C_\nu^\tau d\log(a_\nu - a_\tau)$. After possible introduction of some signs, the $d\log$ terms can be pulled to the right. Suppose, among the $\{\tau_1, \tau_2, \ldots, \tau_{r_{2s}-1}\}$ we have $\tau_{j_1}, \ldots, \tau_{j_a} \neq \tau$ and all the other $\tau_k = \tau$. Then that particular summand on the right multiplies

$$d\log(a_\tau - a_{\tau_{j_1}}) \wedge \cdots \wedge d\log(a_{\tau_{j_a}} - a_\tau) = 0.$$

(Note $x_1 + \cdots + x_{a+1} = 0 \Rightarrow dx_1 \wedge \cdots \wedge dx_{a+1} = 0$.) Thus, one term on the right is nonzero, and

$$M(B)_\tau^\tau = M(B_\tau^\tau),$$

proving the lemma. $\qquad\square$

Since $\mathrm{N}w_n(\Psi)$ is a sum of terms $\mathrm{Tr}\ (M(B))$ as in the lemma, we conclude

$$(4.19) \qquad \mathrm{N}w_n(\Psi) = \sum_{\tau=1}^{\delta} \mathrm{N}w_n(\Psi_\tau),$$

where $\Psi_\tau$ is the connection on $F^{\oplus N}$ given (with notation as above) by

$$1_j \mapsto \sum_{k=1}^{N} 1_k \otimes (\eta_j^k + \sum_{\theta \neq \tau} A_j^{k\theta} d\log(a_\tau - a_\theta)).$$

The connection matrix for $\Psi_\tau$ is thus

$$(4.20) \qquad \Phi + \sum_{\theta \neq \tau} A^\theta d\log(a_\tau - a_\theta)$$

where $\Phi = (\eta_j^k)$ and $A^\theta = (A_j^{k\theta})$.

We now consider the right-hand side of formula $(0.15)$, which in our case takes the form

$$-\mathrm{N}w_n(E, \nabla) \cdot c_1(\Omega_{P/F}^1(\log D), \mathrm{res}_D).$$

Since $\Omega_{P/F}^1(\log D)$ has rank 1, the relative Chern class can be computed in a standard way to be the divisor of any meromorphic section $\omega$ of the bundle such that $\omega$ is regular along $D$ and $\mathrm{res}_D(\omega) = 1$. We shall assume that $0 \notin D$. (This is easy to arrange by application of an automorphism to $P$.) We take for our meromorphic section

$$(4.21) \qquad \omega := \Big( \sum_{\tau=1}^{\delta} \frac{1}{z - a_\tau} - \frac{\delta+1}{z} \Big) dz.$$

Clearing denominators, we get

$$(4.22) \qquad \omega = \frac{F(z)}{z \prod(z - a_\tau)} dz; \quad F(z) = \sum_{\tau=1}^{\delta} a_\tau \prod_{\theta \neq \tau} (z - a_\theta) - \prod_{\tau=1}^{\delta} (z - a_\tau).$$

Formally, with $F = \prod_{i=1}^{\delta}(z - \beta_i)$,

$$c_1(\Omega_{P/F}^1(\log D), \mathrm{res}) = (\omega) = \sum (\beta_i) - (0).$$

We shall need to compute

$$(4.23) \qquad -\sum_{i=1}^{\delta} \mathrm{N}w_n(E, \nabla)|_{z=\beta_i} + \mathrm{N}w_n(E, \nabla)|_{z=0}$$

and compare the answer to $\mathrm{N}w_p(\Phi) - \mathrm{N}w_p(\Psi)$ (cf. $(4.13)$, $(4.15)$, $(4.17)$, $(4.19)$, $(4.20)$).

Now we show some identities which will be used to transform the right-hand side of the Riemann-Roch formula (0.1).

**PROPOSITION 4.4.** *Suppose* $\{j_1, \ldots, j_r\} \subset \{1, \ldots, \delta\}$. *Then*

$$d\log(z - a_{j_1}) \wedge \cdots \wedge d\log(z - a_{j_r})|(\omega)$$

$$:= \sum_{i=1}^{\delta} d\log(\beta_i - a_{j_1}) \wedge \cdots \wedge d\log(\beta_i - a_{j_r}) - d\log(a_{j_1}) \wedge \cdots \wedge d\log(a_{j_r})$$

$$= \sum_{t \neq j_1, \ldots, j_r} d\log(a_t - a_{j_1}) \wedge \cdots \wedge d\log(a_t - a_{j_r}).$$

*Proof.* To simplify notation we consider the case $j_k = k, \ 1 \leq k \leq r$.

**LEMMA 4.5.**

$$\sum_{i=1}^{\delta} d\log(\beta_i - a_1) \wedge \cdots \wedge d\log(\beta_i - a_r)$$

$$= \sum_{j=1}^{r} (-1)^{j-1} d\log(F(a_j)) \wedge d\log(a_j - a_1) \wedge \cdots \wedge d\widehat{\log(a_j - a_j)} \wedge \cdots$$

$$\cdots \wedge d\log(a_j - a_r).$$

*Proof.* We can argue universally. Consider the rings

$$A = \mathbb{Q}[z_1, \ldots, z_\delta, t_1, \ldots, t_r] \subset B \ = \ A[x]/(x^\delta + \sum_{i=1}^{\delta} z_i x^{\delta-i})$$

$$\cong \ \mathbb{Q}[z_1, \ldots, z_{\delta-1}, t_1, \ldots, t_r, X].$$

Note that both $A$ and $B$ are polynomial rings over $\mathbb{Q}$. Let $L \subset M$ be their quotient fields, and consider the symbol

$$S = \{X - t_1, \ldots, X - t_r\} \in K_r(M).$$

We will compute the norm, $N(S) \in K_r(L)$. Write $Z = \operatorname{Spec} B$. The symbol tame : $K_r(M) \to \oplus_{x \in Z^{(1)}} K_{r-1}(\mathbb{Q}(x))$, when applied to $S$ yields

$$\operatorname{tame}(S) = \sum_{k=1}^{r} (-1)^{k-1} \{t_k - t_1, \ldots, \widehat{t_k - t_k}, \ldots, t_k - t_r\}|_{X=t_k}.$$

Let $\pi : \operatorname{Spec}(B) \to \operatorname{Spec}(A)$. We have a map on divisors

$$\pi_*(X - t_i = 0) = (F(t_i) = 0)$$

with degree 1, so that

$$
\begin{aligned}
\mathrm{tame}(N(S)) &= N(\mathrm{tame}(S)) \\
&= \sum_{k=1}^{r}(-1)^{k-1}\{t_k - t_1, \ldots, \widehat{t_k - t_k}, \ldots, t_k - t_r\}|_{F(t_k)=0} \\
&= \mathrm{tame}\Big(\sum_{k=1}^{r}(-1)^{k-1}\{F(t_k), t_k - t_1, \ldots, \widehat{t_k - t_k}, \ldots, t_k - t_r\}\Big).
\end{aligned}
$$

(The last equality holds because $F(t_j)/F(t_k) = 1$ on the divisor $t_j = t_k$.). Since $L$ is purely transcendental over $\mathbb{Q}$, this determines $N(S)$ up to constant symbols, which can be ignored because we want to apply $d\log$. Specializing the $z_i$ to the coefficients of our $F$ and the $t_i \mapsto a_i$ and applying $d\log$, we deduce the lemma. $\qquad\square$

LEMMA 4.6.

$$
\begin{aligned}
d\log(b_1) \wedge \cdots \wedge d\log(b_r) &= \sum_{k=1}^{\delta}(-1)^{k-1}d\log(b_k) \wedge d\log(b_k - b_1) \wedge \cdots \\
&\quad \cdots \wedge d\log\widehat{(b_k - b_k)} \wedge \cdots \wedge d\log(b_k - b_r).
\end{aligned}
$$

*Proof.* As above, we argue universally and prove the corresponding identity for symbols. For this it suffices to compare the images under the tame symbol. At the divisor $b_j - b_k = 0$ for $j < k$ we need

$$
\begin{aligned}
0 &= (-1)^{j+k}\{b_k, b_k - b_1, \ldots, \widehat{b_k - b_j}, \ldots, \widehat{b_k - b_k}, \ldots, b_k - b_r\}|_{b_k = b_j} \\
&\quad + (-1)^{j+k-1}\{b_j, b_j - b_1, \ldots, \widehat{b_j - b_k}, \ldots, \widehat{b_j - b_j}, \ldots, b_j - b_r\}|_{b_k = b_j},
\end{aligned}
$$

which is clear. Finally at the divisor $b_k = 0$ we need

$$
(-1)^{k-1}\{-b_1, \ldots, \widehat{-b_k}, \ldots, -b_r\} = (-1)^{k-1}\{b_1, \ldots, \hat{b}_k, \ldots, b_r\} + \varepsilon
$$

where $\varepsilon$ dies under $d\log$. Again this is clear. $\qquad\square$

Returning to the proof of Proposition 4.4, we apply Lemmas 4.5 and 4.6 (with $b_j = a_j$) to conclude

$$
\begin{aligned}
d\log(z - a_{j_1}) &\wedge \cdots \wedge d\log(z - a_{j_r})|_{(\omega)} \\
&= \sum_{s=1}^{r}(-1)^{s-1}d\log\Big(\prod_{k \notin \{j_1, \ldots, j_r\}}(a_{j_s} - a_k)\Big) \wedge d\log(a_{j_s} - a_{j_1}) \wedge \cdots \\
&\quad \cdots \wedge d\log\widehat{(a_{j_s} - a_{j_s})} \wedge \cdots \wedge d\log(a_{j_s} - a_{j_r}) \\
&= \sum_{\substack{k \neq j_1, \ldots, j_r \\ s=1}}^{s=r}(-1)^{s-1}d\log(a_{j_s} - a_k) \wedge d\log(a_{j_s} - a_{j_1}) \wedge \cdots \\
&\quad \cdots \wedge d\log\widehat{(a_{j_s} - a_{j_s})} \wedge \cdots \wedge d\log(a_{j_s} - a_{j_r}).
\end{aligned}
$$

Finally we apply Lemma 4.5 again to this last expression, taking $b_s = a_{j_s} - a_k$, to get the assertion of the proposition:

$$d\log(z-a_{j_1})\wedge\cdots\wedge d\log(z-a_{j_r})|_{(\omega)} = \sum_{k\neq j_1,\ldots,j_r} d\log(a_{j_1}-a_k)\wedge\cdots\wedge d\log(a_{j_r}-a_k).$$

PROPOSITION 4.7. *With notation as above, formula* (0.15) *holds for* $(E, \nabla)$.

*Proof.* The computation mentioned in (4.23) can be done as follows. Let $\rho_\nu$ be closed 1-forms. For $J = \{j_1 < \ldots < j_r\} \subset \{1,\ldots,\delta\}$ define $\rho_J = \rho_{j_1} \wedge \cdots \wedge \rho_{j_r}$. Write

(4.24)
$$\mathrm{Nw}_n(\mathcal{O}_P^{\oplus N}, \sum_{\nu=1}^{\delta} A^\nu \rho_\nu + \Phi) = \sum_{J\subset\{1,\ldots,\delta\}} P_J(A^\nu, dA^\nu, \Phi, d\Phi)\rho_J + \mathrm{Nw}_n(\mathcal{O}_P^{\oplus N}, \Phi).$$

Here $A^\nu$ (resp. $\Phi$) are matrices with coefficients in $F$ (resp. $\Omega_F^1$), and the $P_J$ are independent of the $\rho_j$. Then, using Proposition 4.4, we get

(4.25)  $-\mathrm{Nw}_n(\sum_{\nu=1}^{\delta} A^\nu \rho_\nu + \Phi)|_{(\omega)}$

$$= -\sum_{\substack{J\subset\{1,\ldots,\delta\}\\r=|J|\geq1}} P_J(A^\nu, dA^\nu, \Phi, d\Phi) \sum_{k\notin J} d\log(a_{j_1} - a_k) \wedge \cdots$$

$$\cdots \wedge d\log(a_{j_r} - a_k) + (1-\delta)\mathrm{Nw}_n(F^{\oplus N}, \Phi).$$

On the other hand, if we fix $\tau \leq \delta$ and take $\rho_\nu = d\log(a_\tau - a_\nu)$ for $\nu \neq \tau$ and $\rho_\tau = 0$ we find

(4.26)  $\mathrm{Nw}_n(\mathbb{R}f_*(E \otimes \Omega_{P/F}^*(\log D))) = \mathrm{Nw}_n(F^N, \Phi) - \sum_{\tau=1}^{\delta} \mathrm{Nw}_n(F^N, \Psi_\tau)$

$$= -\sum_{\tau=1}^{\delta} \sum_{\substack{J\subset\{1,\ldots,\delta\}\\\tau\notin J}} P_J(A^\nu, dA^\nu, \Phi, d\Phi) d\log(a_{j_1} - a_\tau) \wedge \cdots$$

$$\cdots \wedge d\log(a_{j_r} - a_\tau) + (1-\delta)\mathrm{Nw}_n(\Phi).$$

The right-hand sides of (4.25) and (4.26) coincide, proving the proposition. □

We have assumed throughout that the divisor $D$ is a sum of $F$-rational points. In proving the Riemann-Roch theorem for $\mathrm{Nw}_n$ with $n \geq 2$ this is not a problem. These classes take values in a group without torsion. We may argue as in the proof of Reduction 3.7 and pull back to a finite field extension $F'/F$. For $\mathrm{Nw}_1$ we must be more careful.

LEMMA 4.8. *Let $L/F$ be a finite, Galois extension of fields, with Galois group $G$. Let $S$ be a finite $G$-set, and let $L[S]$ be the $L$-vector space spanned by $s \in S$. Then the natural map*

$$(L[S])^G \otimes_F L \to L[S]$$

*is an isomorphism, where $G$ acts on $L[S]$ by $g(\sum \ell_s[s]) = \sum g(\ell_s)[g(s)]$.*

*Proof.* The normal basis theorem gives $L \cong F[G]$ as a $G$-module, so that $L[S] \cong F[G \times S]$ with $G$ acting diagonally. The set $e \times S \subset G \times S$ is a set of coset representatives for the diagonal action of $G$ on $G \times S$, so that $(L[S])^G$ has $F$-dimension $|S|$, and it suffices to show the above map is injective. Let $x_1, \ldots, x_\delta$ be an $F$-basis for $(L[S])^G$, and let $\sum \ell_i x_i \mapsto 0$ be a nonzero element in the kernel with the minimal number of nonzero $\ell_i$. We may assume $\ell_1 = 1$. If some $\ell_i \notin F$ we can find $g$ such that $g(\ell_i) \neq \ell_i$ and observe that $\sum(g(\ell_i) - \ell_i)x_i$ is a nontrivial element in the kernel with fewer nonzero $\ell_i$. $\square$

In proving the Riemann-Roch theorem, we have reduced to the case $E \cong \mathcal{O}_P^{\oplus N}$, and, as remarked at the beginning of this section, this isomorphism can be taken to be defined over $F$. Let $L$ be the Galois extension of $F$ generated by the $a_\nu$, $1 \leq \nu \leq \delta$, and write $G$ for the Galois group. From the connection equation (4.8) it follows that $\Phi = (\eta_j^k)$ is a matrix with coefficients in $\Omega_F^1$. Applying the lemma to the set $\{d\log(z - a_\nu)\}$, $1 \leq \nu \leq \delta$, we find an invertible $\delta \times \delta$-matrix $(\alpha_\ell^\nu)$ in $L$ such that the logarithmic forms

$$(4.27) \qquad \sum_{\nu=1}^{\delta} \alpha_\ell^\nu d\log(z - a_\nu), \quad 1 \leq \ell \leq \delta,$$

form an $F$-basis for $\Gamma(P, \Omega_{P/F}^1(\log D))$. With respect to this new basis, the connection matrix $B$ for $\Psi$ (4.15) becomes

$$(4.28) \qquad \beta B \beta^{-1} + d\beta \beta^{-1}$$

where $\beta := (\alpha_\ell^\nu) \otimes I_N$, with $I_N$ the $N \times N$ identity matrix. We conclude

$$(4.29) \quad w_1(\mathbb{R}f_* E \otimes \Omega_{P/F}^*(\log D)) \quad = \quad -w_1(E, \nabla) \cdot c_1(\Omega_{P/F}^1(\log D), \mathrm{res})$$
$$+ N \cdot d\log(\det(\alpha_\ell^\nu)).$$

Replacing $E$ with the rank 0 virtual bundle $E - N \cdot \mathcal{O}_P$, we get the desired Riemann-Roch theorem in this case:

$$(4.30) \qquad w_1(\mathbb{R}f_*(E - N \cdot \mathcal{O}) \otimes \Omega_{P/F}^*(\log D))$$
$$= -w_1(E - N \cdot \mathcal{O}, \nabla - N \cdot d) \cdot c_1(\Omega_{P/F}^1(\log D), \mathrm{res}).$$

It follows from $G$-invariance of the form in (4.27) that

$$(4.31) \qquad g(\alpha_\ell^\nu) = \alpha_\ell^{g(\nu)} \text{ where we define } g(\nu) \text{ by } a_{g(\nu)} = g(a_\nu).$$

To verify the desired 2-torsion condition, we remark that the matrix $\alpha \cdot {}^t\alpha$ has entries $\sum_\nu \alpha_\ell^\nu \alpha_m^\nu \in F$ by (4.31), so

$$(4.32) \qquad 2d \log(\det \alpha) = 0 \text{ in } \Omega_F^1/F^\times = H_{\mathrm{CS}}^2(\mathrm{Spec}(F)).$$

## 5. Connection on the determinant line for curves

It is curious that in the basic case of a curve over a function field, the Riemann-Roch theorem for $\mathrm{N}w_1 = w_1$ does not require the basic curvature form condition.

THEOREM 5.1. *Let $f : X \to S = \mathrm{Spec}(F)$ be a smooth, complete curve over a function field. Let $D \subset X$ be a reduced, effective divisor, and let*

$$\nabla : E \to E \otimes \Omega_X^1(\log D)$$

*be a vector bundle of rank $N$ with connection. Then there is a naturally defined connection on the determinant bundle*

$$\det(\mathbb{R}f_*((E - N \cdot \mathcal{O}) \otimes \Omega_{X/S}^*(\log D)))$$

*and the Riemann-Roch formula holds for line bundles with connection:*

$$\det(\mathbb{R}f_*((E - N \cdot \mathcal{O}) \otimes \Omega_{X/S}^*(\log D)) = f_*\bigg(\det(E) \cdot c_1(\Omega_{X/S}^1(\log D), \mathrm{res}_D)\bigg).$$

*Proof.* We assume for a while that

$$(5.1) \qquad R^1 f_* E = R^1 f_*(\Omega_{X/S}^1(\log D) \otimes E) = 0.$$

Then the complex of $F$-vectorspaces

$$f_* E \to f_*(\Omega_{X/S}^1(\log D) \otimes E)$$

represents

$$Rf_*(\Omega_{X/S}^*(\log D) \otimes E).$$

Let

$$\sigma : f_*(\Omega_{X/S}^1(\log D) \otimes E) \to f_*(\Omega_X^1(\log D) \otimes E)$$

be a splitting of the exact sequence

$$0 \to \Omega_S^1 \otimes f_* E \to f_*(\Omega_X^1(\log D) \otimes E) \to f_*(\Omega_{X/S}^1(\log D) \otimes E) \to 0.$$

This gives rise to the diagram 4.10 with $P/F$ replaced by $X/S$, $\Phi = \tau \circ \nabla$, $\Psi = \nabla_{X/S} \circ \sigma$, except that in our situation, $\nabla_2 \circ \nabla_1 \neq 0$ if the curvature does not fulfill Definition 3.1. This defines the following diagram, similar to

4.12, except that it does not commute if the curvature condition (3.1) is not satisfied:

$$
\begin{array}{ccc}
\Gamma(E) & \xrightarrow{\ \Phi\ } & \Gamma(E) \otimes \Omega^1_S \\[2pt]
\ \ \downarrow{\scriptstyle \nabla_{X/S}} & & \ \ \downarrow{\scriptstyle \nabla_{X/S}\otimes 1} \\[4pt]
\Gamma(E \otimes \Omega^1_{X/S}(\log D)) & \xrightarrow{\ \Psi\ } & \Gamma(E \otimes \Omega^1_{X/S}(\log D)) \otimes \Omega^1_S
\end{array}
$$

(5.2)

with

$$
(5.3) \qquad \Phi := \nabla_1 - \sigma \circ \nabla_{X/F}, \Psi = -\nabla_2 \circ \sigma.
$$

PROPOSITION 5.2. *The connection* $\mathrm{GM}(\nabla) = w_1(\Phi) - w_1(\Psi)$ *on*

$$
\det \left( Rf_*(\Omega^*_{X/S}(\log D) \otimes E) \right)
$$

*is well-defined, and Theorem* 5.1 *holds true for all coherent sheaves with connections* $(E, \nabla)$ *fulfilling condition* (5.1).

*Proof.* Another splitting is of the shape $\sigma' = \sigma + \varphi$, where

$$
\varphi : \Gamma(\Omega^1_{X/S}(\log D) \otimes E) \to \Omega^1_S \otimes \Gamma(E)
$$

is an $F$-linear map. Thus

$$
(w_1(\Phi') - w_1(\Psi')) - (w_1(\Phi) - w_1(\Psi)) = \mathrm{Tr}(\varphi \circ \nabla_{X/S} - \nabla_{X/S} \circ \varphi \otimes 1) = 0.
$$

Once the connection on the determinant line is defined, one has to verify that one can apply Reduction 3.5 and Section 4. Take $Y = D + H$ containing the ramification of $g$, with $H \cap D = \emptyset$, where $g : X \to \mathbb{P}^1_S$ is as in Reduction 3.5. Then of course $R^1 f_*(\Omega^1_{X/S}(\log Y) \otimes E) = 0$. On the other hand, condition (5.1) implies that $R^1 f_*(\Omega^1_X(\log D) \otimes E) = 0$. Thus one has a commutative diagram of exact sequences

(5.4)

$$
\begin{array}{ccccccccc}
0 & \to & f_*(\Omega^1_X(\log D) \otimes E) & \to & f_*(\Omega^1_X(\log Y) \otimes E) & \to & f_*E|H & \to & 0 \\[2pt]
& & \downarrow & & \downarrow & & \ \|\ \downarrow & & \\[4pt]
0 & \to & f_*(\Omega^1_{X/S}(\log D) \otimes E) & \to & f_*(\Omega^1_{X/S}(\log Y) \otimes E) & \to & f_*E|H & \to & 0.
\end{array}
$$

One chooses a splitting

$$
\sigma' : f_*(\Omega^1_{X/S}(\log Y) \otimes E) \to f_*(\Omega^1_X(\log Y) \otimes E)
$$

with

$$
(5.5) \qquad \sigma'|f_*(\Omega^1_{X/S}(\log D) \otimes E) = \sigma,
$$

$$
\sigma' \bmod \sigma = \mathrm{Id} : f_*E|H \to f_*E|H.
$$

This induces $\Phi' = \Phi, \Psi'$, and condition (5.5) implies that one has an exact sequence of $F$-connections

$$0 \to \Psi \to \Psi' \to f_*\nabla|H \to 0.$$

Thus one obtains

$$(5.6) \qquad Nw_1(\Phi) - Nw_1(\Psi) = Nw_1(\Phi) - Nw_1(\Psi') + Nw_1(f_*\nabla|H).$$

This shows that Theorem 5.1 for $(E, \nabla, D)$ is equivalent to Theorem 5.1 for $(E, \nabla, Y)$. Now we can apply Reduction 3.5. Moreover, since in Section 4, $GM^i(\nabla)$ was described via the diagrams (4.10) and (4.12), this concludes the proof of the proposition. $\qquad \square$

Let $(E, \nabla)$ be any connection on $X$ as in Theorem 5.1. Let $Y = D + H$, with $H \cap D = \emptyset$ such that condition (5.1) is fulfilled with $E$ replaced by $E(H)$ and $D$ replaced by $Y$. Then, by [9] (for an algebraic version of it, see e.g. [15]) the inclusion

$$\Omega^*_{X/S}(\log Y) \otimes E \to \Omega^*_{X/S}(\log Y) \otimes E(H)$$

is a quasi-isomorphism. We may thus define

$$Nw_1(Rf_*\Omega^*_{X/S}(\log Y) \otimes E) := Nw_1(f_*(\Omega^*_{X/S}(\log Y) \otimes E(H))). \qquad \square$$

PROPOSITION 5.3. *The class*

$$Nw_1(Rf_*\Omega^1_{X/S}(\log Y) \otimes E) + Nw_1(f_*(E|H), f_*(\nabla|H))$$

*does not depend on the choice of $H$.*

*Proof.* Since the condition (5.1) for $E(H)$ and $Y = D + H$ implies the condition (5.1) for $E(H + K)$ and $Y + K$ for any effective divisor $K$, it is sufficient to show

$$(5.7) \qquad Nw_1(Rf_*\Omega^1_{X/S}(\log Y) \otimes E) + Nw_1(f_*\nabla|H)$$
$$= Nw_1(Rf_*\Omega^1_{X/S}(\log(Y + K)) \otimes E) + Nw_1(f_*\nabla|(H + K))$$

for either an irreducible component $K$ of $H$ or for an irreducible divisor $K$ disjoint of $Y$ to show that the class is well-defined. The first case is trivial and the second case is treated as the proof of Proposition 5.2. Theorem 5.1 now follows from Proposition 5.2. $\qquad \square$

UNIVERSITY OF CHICAGO, CHICAGO IL
*E-mail address*: bloch@math.uchicago.edu

UNIVERSITÄT ESSEN, ESSEN, GERMANY
*E-mail address*: esnault@uni-essen.de